\documentclass[]{article}
\usepackage[left=2.5cm, right=2.5cm, top=2cm]{geometry}
\usepackage{blkarray}
\usepackage{dsfont}
\usepackage{bm}
\usepackage[cmex10]{amsmath}
\usepackage{amsmath}
\usepackage{amssymb}
\usepackage{epsfig,floatflt}
\usepackage{psfrag,graphicx}
\usepackage[colorlinks]{hyperref}
\usepackage{cite}
\usepackage{epstopdf}
\usepackage{graphicx}
\usepackage{subcaption}

\newcommand{\bx}{\bm{\bar{x}}}
\newcommand{\bu}{\bm{\bar{u}}}
\newcommand{\blam}{\bm{\bar{\lambda}}}
\newcommand{\bmu}{\bm{\bar{\mu}}}
\newcommand{\IR}{\mathbb{R}}

\DeclareMathOperator*{\argmin}{arg\,min}

\newtheorem{thm}{\textbf{Theorem}}
\newtheorem{prob}{\textbf{Problem}}
\newtheorem{prop}{\textbf{Property}}
\newtheorem{lem}{\textbf{Lemma}}

\newtheorem{cor}{\textbf{Corollary}}
\newtheorem{asm}{\textbf{Assumption}}
\newtheorem{rem}{\textbf{Remark}}
\newtheorem{exm}{\textbf{Example}}

\DeclareSymbolFont{extraup}{U}{zavm}{m}{n}
\DeclareMathSymbol{\vardiamond}{\mathalpha}{extraup}{87}
\newcommand*{\qedlem}{\hfill\ensuremath{\blacksquare}}%
\newcommand*{\qedprob}{\hfill\ensuremath{\triangledown}}%
\newcommand*{\qedprop}{\hfill\ensuremath{\square}}%
\newcommand*{\qedproof}{\hfill\ensuremath{\spadesuit}}%

\begin{document}

\title{Bernstein approximation of optimal control problems\thanks{This work was supported by AFOSR, ONR, NSF and NASA.}}

\author{Venanzio Cichella\thanks{Venanzio Cichella is with the Department of Mechanical Engineering, University of Iowa, Iowa City, 52242 IA {\tt\small venanzio-cichella@uiowa.edu}}, 
Isaac Kaminer\thanks{Isaac Kaminer and Claire Walton are with the Department of Mechanical and Aerospace Engineering, Naval Postgraduate School, Monterey, CA 93940 {\tt\small \{kaminer, clwalton1 \}@nps.edu}}, 
Claire Walton\footnotemark[3], 
Naira Hovakimyan\thanks{Naira Hovakimyan is with the Department of Mechanical Science and Engineering, University of Illinois at Urbana-Champaign, Urbana, IL 61801 {\tt\small nhovakim@illinois.edu}}, Ant{\'{o}}nio M. Pascoal\thanks{Ant{\'{o}}nio M. Pascoal is with the Institute for Systems and Robotics (ISR), Instituto Superior Tecnico (IST), Univ. Lisbon, Portugal. {\tt\small antoniog@isr.ist.utl.pt}}%
}
%




\maketitle

\begin{abstract}
Bernstein polynomial approximation to a continuous function has a slower rate of convergence as compared to other approximation methods. \emph{``The fact seems to have precluded any numerical application of Bernstein polynomials from having been made. Perhaps they
will find application when the properties of the approximant in the large are of more importance than the closeness of the approximation.''} -- has remarked P.J. Davis in his 1963 book Interpolation and Approximation.

This paper presents a direct approximation method for nonlinear optimal control problems with mixed input and state constraints based on Bernstein polynomial approximation. We provide a rigorous analysis showing that the proposed method yields consistent approximations of time continuous optimal control problems. Furthermore, we demonstrate that the proposed method can also be used for costate estimation of the optimal control problems. This latter result leads to the formulation of the Covector Mapping Theorem for Bernstein polynomial approximation. Finally, we explore the numerical and geometric properties of Bernstein polynomials, and illustrate the advantages of the proposed approximation method through several numerical examples.
\end{abstract}

\section{Introduction}
\label{sec:Introduction}
Optimal control problems that arise from most engineering applications are in general very complex. Finding a closed-form solution to these problems can be difficult or even impossible, and therefore they must be solved numerically.
Numerical methods include indirect and direct methods \cite{rao2009survey}. Indirect methods solve the problems by converting them into boundary value problems. Then, the solutions are found by solving systems of differential equations. On the other hand, direct methods are based on transcribing optimal control problems into nonlinear programming problems (NLPs) using some discretization scheme \cite{rao2009survey,betts2010practical,betts1998survey,conway2012survey}. These NLPs can be solved using ready-to-use NLP solvers (e.g. MATLAB, SNOPT, etc.) and do not require calculation of costate and adjoint variables as indirect methods do.

A pioneering work in the literature on direct methods is the one of Polak on consistency of approximation theory reported in his book (see \cite[Section 3.3]{elijah1997optimization}). Borrowing tools from variational analysis, Polak provides a theoretical framework to assess the convergence properties of direct methods. Motivated by the consistency of approximation theory, a wide range of direct methods that use different discretization schemes have been developed, including Euler \cite{elijah1997optimization}, Runge-Kutta \cite{schwartz1996consistent}, Pseudospectral \cite{ross2012review} methods, as well as the method presented in this paper.

Pseudospectral methods are the most popular direct methods nowadays, and they have been applied successfully for solving a wide range of optimization problems, e.g. \cite{fahroo2006discrete,bollino2007pseudospectral,gong2009pseudospectral,bedrossian2009zero,bollino2008collision,bedrossian2007first,ross2012review}. They offer several advantages over many other discretization methods, mainly owing to their spectral (exponential) rate of convergence. However, as pointed out in \cite{chen2015decoupled,augugliaro2012generation,ChoePhd,cichella2018optimal}, there is one salient disadvantage associated with these methods. When discretizing the state and/or the input, the constraints are enforced at the discretization nodes; unfortunately, satisfaction of constraints cannot be guaranteed in between the nodes. To avoid violation of the constraints in between the nodes, the order of approximation (number of nodes) can be increased; however, this leads to larger NLPs, which may become computationally expensive and too inefficient to solve. This problem does not limit itself to pseudospectral methods, but it is common to methods that are based on discretization.

This undesirable behaviour becomes obvious, for example, when considering the optimal trajectory generation problem for multi-vehicle missions, where a large number of vehicles have to reach their final destinations by following trajectories that have to guarantee intervehicle separation for safety all the time. Clearly, with a small order of approximation, separation between the trajectories will be hardly satisfied. Increasing the number of nodes will eventually produce spatially separated trajectories, but will also drastically increase the number of collision avoidance constraints and thus also the complexity of the problem.

Pseudospectral methods also suffer from a drawback when dealing with non-smooth optimal control problems. This drawback is mainly related to the well known Gibbs phenomenon \cite{hewitt1979gibbs}, common to all approximation methods based on orthogonal polynomials. The Gibbs phenomenon, visible in the form of oscillations, reduces the accuracy of the approximation to first order away from discontinuities and to $\mathcal{O}(1)$ in the neighborhoods of jumps \cite{gelb2006robust}.
Several extensions of pseudospectral methods have been developed to deal with this disadvantage and lessen the effect of the Gibbs phenomenon (e.g. \cite{tohidi2013efficient,ross2004pseudospectral,darby2011hp}). The most accurate methods require the location of the discontinuities to be known a priori, which is often impractical or difficult. Other methods are based on the estimation of these locations, which could result in inefficiency or ill conditioning of the discretized problem, especially when the number of discontinuities is large and unknown.

The present article proposes a direct method based on Bernstein approximation. Bernstein approximants have several nice properties. First of all, the approximants converge uniformly to the functions that they approximate -- and so do their derivatives \cite[Chapter 3]{farouki2012bernstein}. This, as we will discuss later, is useful for derivation of convergence properties of the proposed computational method. Moreover, Bernstein polynomials behave well, even when the functions being approximated are non-smooth. In fact, as demonstrated in \cite{gzyl2003approximation}, the Gibbs phenomenon does not occur when approximating piecewise smooth, monotone functions with both left and right derivatives at every point by Bernstein polynomials. As a result, the proposed method based on Bernstein approximation lends itself to problems that have discontinuous states and/or controls, e.g. bang-bang optimal control problems (see also \cite{ricciardi2018direct}).
Finally, due to their favorable geometric properties (see \cite[Chapter 5]{farouki2012bernstein}) Bernstein polynomials afford computationally efficient algorithms for the computation of state and input constraints for the whole time interval where optimization takes place, and not only at discretization points (see \cite{chang2011computation,ChoePhd}). Hence, with the proposed approach the solutions can be guaranteed to be feasible and satisfy the constraints for all times, while retaining the computational efficiency of methods based on discretization.

\emph{``There is a price that must be paid for these beautiful approximation properties: the convergence of the Bernstein polynomials is very slow.''} -- 
wrote P.J. Davis in his 1963 book Interpolation and Approximation \cite{davis1963interpolation}. He continues:
\emph{``This fact seems to have precluded any numerical application of Bernstein polynomials from having been made. Perhaps they will find application when the properties of the approximant in the large are of more importance than the closeness of the approximation.''}
In fact, the slow convergence of the Bernstein approximation implies that the approach proposed in the present paper is outperformed by, for example, pseudospectral methods in terms of convergence rate. This is not surprising, since the choice of nodes and the interpolating polynomials in the pseudospectral methods is dictated by approximation accuracy and convergence speed, while sacrificing satisfaction of constraints in between the nodes. On the other hand, our approach prioritizes constraint satisfaction at the expense of a slower convergence rate.

The paper is structured as follows: in Section~\ref{sec:MathematicalPreliminaries} we present the notation and the mathematical results, which will be used later in the paper. Section~\ref{sec:ProblemFormulation} introduces the optimal control problem of interest and some related assumptions.
Section~\ref{sec:BernApprox} presents the NLP method based on Bernstein approximation that approximates the optimal control problem. Section~\ref{sec:analysis} demonstrates that the proposed approximation method yields approximate results that converges uniformly to the optimal solution. In Section \ref{sec:costate} we derive the Karush–Kuhn–Tucker (KKT) conditions associated with the NLP. Section \ref{sec:costateresults} compares these conditions to the first order optimality conditions for the original optimal control problem and states the Covector Mapping Theorem for Bernstein approximation. Numerical examples are discussed in Section~\ref{sec:numericalresults}. The paper ends with conclusions in Section~\ref{sec:conclusions}.


\section{Notation and mathematical background}\label{sec:MathematicalPreliminaries}
Vector valued functions are denoted by bold letters, $\bm{x}(t) =[x_1(t) \, , \, \ldots \, , \, x_n(t)]^\top$, while vectors are denoted by bold letters with an upper bar, $\bx = [x_1 \, , \, \ldots \, , \, x_n]^\top \in \IR^{n}$. The symbol $\mathcal{C}^r$ denotes the space of functions with $r$ continuous derivatives. $\mathcal{C}^r_n$ denotes the space of $n$-vector valued functions in $\mathcal{C}^r$. $||\cdot||$ denotes the Euclidean norm, $||\bx|| = \sqrt{x_1^2 + \ldots + x_n^2}$.

\subsection{Bernstein polynomials}

The Bernstein basis polynomials of degree $N$ are defined as
$$ b_{j,N}(t) = \binom{N}{j} t^j(1-t)^{N-j}  \, , \qquad t \in [0,1] \, ,$$
for $j=0,\ldots,N$, with
$$
\binom{N}{j} = \frac{N !}{j!(N-j)!}\,.
$$
They were originally introduced by the mathematician Sergei Natanovich Bernstein in 1912 to facilitate a constructive proof of the Weierstrass approximation theorem \cite{bernstein1912}. An $N$th order Bernstein polynomial $x_N:[0,1]\to\IR$ is a linear combination of $N+1$ Bernstein basis polynomials of order $N$, i.e.
$$
x_N(t) = \sum_{j=0}^N \bar{x}_j b_{j,N}(t) \, , \qquad t \in [0,1] \, ,
$$
where $\bar{x}_j \in \IR$, $j=0,\ldots,N$, are referred to as Bernstein coefficients (also known as control points). For the sake of generality, and with a slight abuse of terminology, in this paper we extend the definition of a Bernstein polynomial given above to a vector of $N$th order polynomials $\bm{x}_N : [0,1] \to \IR^n$ expressed in the following form
\begin{equation} \label{eq:beziercurve}
\bm{x}_N(t) = \sum_{j=0}^N  \bx_{j,N} b_{j,N}(t)  \, , \qquad t \in [0,1] \, ,
\end{equation}
where $\bx_{0,N},\ldots,\bx_{N,N}$ are $n$-dimensional Bernstein coefficients.

Bernstein polynomials were popularized by Pierre B\'{e}zier in the early 1960s as useful tools for geometric design (B\'{e}zier used Bernstein polynomials to design the shape of the cars at the Renault company in France), and are now widely used in computer graphics, animations and type fonts such as postscript fonts and true type fonts. For this reason, the Bernstein polynomial introduced in Equation \eqref{eq:beziercurve} is often referred to as a B{\' e}zier curve, especially when used to describe a spatial curve.

Bernstein polynomials possess favorable geometric and numerical properties that can be exploited in many application domains.
For an extensive review on Bernstein polynomials and their properties the reader is referred to \cite{farouki2012bernstein}.
The derivative and integral of a Bernstein polynomial $\bm{x}_N(t)$ can be easily computed as
\begin{equation*}
\dot{\bm{x}}_N(t) = N \sum_{j = 0}^{N-1} (\bx_{j+1,N}-\bx_{j,N}) b_{j,N-1}(t) \,
\end{equation*}
and
\begin{equation} \label{eq:integralbezier}
\int_0^1{\bm{x}}_N(t) = w \sum_{j = 0}^{N} \bx_{j,N} \, , \qquad w = \frac{1}{N+1} \, ,
\end{equation}
respectively.
Bernstein polynomials can be used to approximate smooth functions.
Consider a $n$-vector valued function $\bm{x}:[0,1] \to \IR^n$. The $N$th order \emph{Bernstein approximation} of $\bm{x}(t)$ is a vector of Bernstein polynomials $\bm{x}_N (t)$ computed as in \eqref{eq:beziercurve} with $\bx_{j,N} = \bm{x}(t_j)$ and $t_j = \frac{j}{N}$ for all $j = 0,\ldots,N$. Namely,
\begin{equation} \label{eq:bernapprox}
\bm{x}_N(t) = \sum_{j=0}^N  \bm{x}(t_j) b_{j,N}(t)  \, , \qquad t_j = \frac{j}{N} \, .
\end{equation}
The following results hold for Bernstein approximations.
\begin{lem}[Uniform convergence of Bernstein approximation]
	\label{lem:bernsteinapprox}
	Let $\bm{x}(t) \in \mathcal{C}^{0}_n$ on $[0,1]$, and $\bm{x}_N(t)$ be computed as in Equation \eqref{eq:bernapprox}. Then, for arbitrary order of approximation $N \in \mathbb{Z}^+$, the Bernstein approximation $\bm{x}_N(t)$ satisfies
	\begin{equation*}
	||\bm{x}_N(t)-\bm{x}(t)|| \leq C_0 W_x(N^{-\frac{1}{2}}) \, ,
	\end{equation*}
	where $C_0$ is a positive constant satisfying $C_0 < 5n/4$, and $W_x(\cdot)$ is the modulus of continuity of $\bm{x}(t)$ in $[0,1]$ \cite{bojanic1989rate,popoviciu1935approximation,sikkema1961wert}. 
	Moreover, if $\bm{x}(t) \in \mathcal{C}^{1}_n$, then
	\begin{equation*}
	\Vert \dot{\bm{x}}_N(t)-\dot{\bm{x}}(t) \Vert \leq C_1 W_{x^\prime}(N^{-\frac{1}{2}}) \, ,
	\end{equation*}
	where $C_1$ is a positive constant satisfying $C_1 < 9n/4$ and $W_{x^{\prime}}(\cdot)$ is the modulus of continuity of $\dot{\bm{x}}(t)$ in $[0,1]$ \cite{powell1981approximation}.
	
\qedlem
\end{lem}
\begin{lem} \cite{floater2005convergence} \label{lem:bernsteinapproxC2}
Assume $\bm{x}(t) \in \mathcal{C}^{r+2}_n$, $r\geq 0$, and let $\bm{x}_N(t)$ be computed as in Equation \eqref{eq:bernapprox}. Let $\bm{x}^{(r)}(t)$ denote the $r$th derivative of $\bm{x}(t)$. Then, the following inequalities hold for all $t\in [0,1]$:
\begin{equation*}
\begin{split}
||\bm{x}_N(t)-\bm{x}(t) || & \leq \frac{C_0}{N} \\
& \, \, \, \vdots \\
||\bm{x}^{(r)}_N(t)-\bm{x}^{(r)}(t) || & \leq \frac{C_r}{N} \, ,
\end{split}
\end{equation*}
 where $C_0,\ldots,C_r$ are independent of $N$.

\qedlem
\end{lem}

\begin{lem} \label{lem:quadrature}
If $\bm{x}(t) \in \mathcal{C}^{0}_n$ on $[0,1]$, then we have
$$
\left\Vert \int_0^1 \bm{x}(t) dt - w \sum_{j=0}^N \bm{x}\left(  \frac{j}{N} \right) \right\Vert \leq C_I W_x(N^{-\frac{1}{2}}) \, 
$$
with $w = \frac{1}{N+1}$, where $C_I> 0 $ is independent of $N$.
Moreover, if $\bm{x}(t) \in \mathcal{C}^{2}_n$, then
$$
\left\Vert \int_0^1 \bm{x}(t) dt - w \sum_{j=0}^N \bm{x}\left(  \frac{j}{N} \right) \right\Vert \leq \frac{C_I}{N} \, .
$$

\qedlem
\end{lem}
The Lemma above follows directly from Lemmas \ref{lem:bernsteinapprox} and \ref{lem:bernsteinapproxC2} and Equation \eqref{eq:integralbezier}.

\begin{figure}
\centering
\begin{subfigure}[b]{0.4\linewidth}
\includegraphics[trim={10cm 0 10cm 0},clip,width=\linewidth]{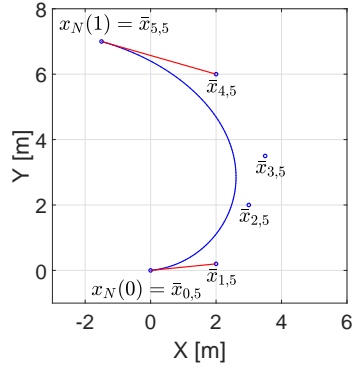}
\caption{End point values property.}\label{app.fig:endpoints}
\end{subfigure}
\begin{subfigure}[b]{0.4\linewidth}
\includegraphics[trim={10cm 0 10cm 0},clip,width=\linewidth]{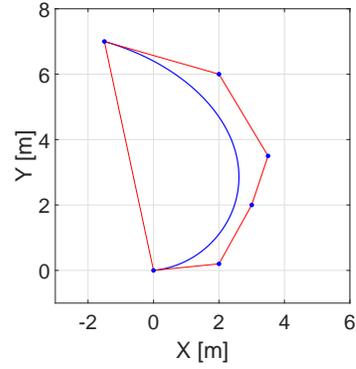}
\caption{Convex hull property.}\label{app.fig:convexhull}
\end{subfigure}

\begin{subfigure}[b]{.4\linewidth}
\includegraphics[trim={10cm 0 10cm 0},clip,width=\linewidth]{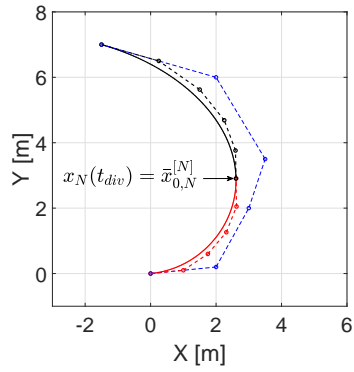}
\caption{de Casteljau algorithm.}\label{app.fig:decasteljau}
\end{subfigure}
\begin{subfigure}[b]{.4\linewidth}
\includegraphics[trim={10cm 0 10cm 0},clip,width=\linewidth]{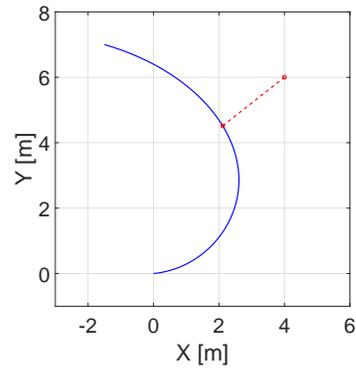}
\caption{Minimum distance to a point.}\label{app.fig:mindistpt}
\end{subfigure}

\begin{subfigure}[b]{.58\linewidth}
\includegraphics[width=\linewidth]{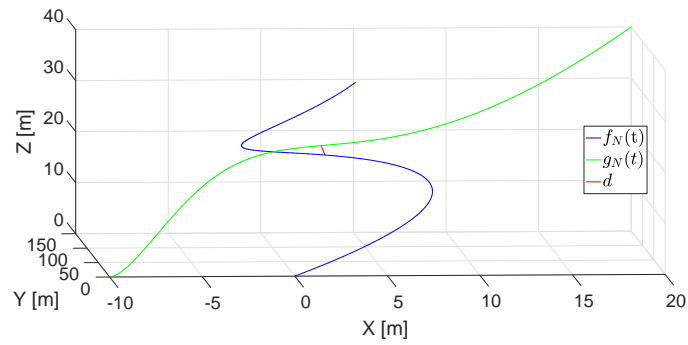}
\caption{Distance between two 3D Bernstein polynomials}\label{app.fig:mindist}
\end{subfigure}
\caption{2D and 3D spatial curves defined by Bernstein polynomials.}
\label{app.fig:bernsteinpoly}
\end{figure}

The following properties of Bernstein polynomials are relevant to this paper.
\begin{prop}[End point values] \label{app.prop:endpoint}
The Bernstein polynomial given by Equation \eqref{eq:beziercurve} satisfies $\bm{x}_N (0) = \bm{\bar{x}}_{0,N}$ and $\bm{x}_N(1) = \bm{\bar{x}}_{N,N}$.
Moreover, the tangent of a Bernstein polynomial at the initial and final points lies on the vectors $\bm{\bar{x}}_{1,N}-\bm{\bar{x}}_{0,N}$ and $\bm{\bar{x}}_{N,N}-\bm{\bar{x}}_{N-1,N}$, respectively. A graphical depiction of this property is illustrated in Figure \ref{app.fig:endpoints}.
\qedprop
\end{prop}

\begin{prop}[Convex hull] \label{app.prop:convexhull}
A Bernstein polynomial is completely contained in the convex hull of its Bernstein coefficients (see Figure \ref{app.fig:convexhull}).

\qedprop
\end{prop}

\begin{prop}[de Casteljau Algorithm] \label{app.prop:decast}
The de Casteljau algorithm is an efficient and numerically stable recursive method to evaluate a Bernstein polynomial at any given point.
The de Casteljau algorithm is also used to split a Bernstein polynomial into two independent ones.
Given an $N$th order Bernstein polynomial $\bm{x}_N: [0,1] \to \IR^d$, and a scalar $t_\emph{div} \in [0,1]$, the Bernstein polynomial at $t_{\emph{div}}$ can be computed using the following recursive relation
$$
\bm{\bar{x}}_{i,N}^{[0]} = \bm{\bar{x}}_{i,N}\, , \quad i = 0,\ldots,N
$$
$$
\bm{\bar{x}}_{i,N}^{[j]} = \bm{\bar{x}}_{i,N}^{[j-1]} (1-t_{\emph{div}})+\bm{\bar{x}}_{i+1,N}^{[j-1]} t_{\emph{div}}  \, , \quad i = 0,\ldots,N-j \, , \quad j = 1,\ldots,N \, .
$$
Then, the Bernstein polynomial evaluated at $t_{\emph{div}}$ is given by
$$
\bm{x}_N(t_{\emph{div}}) = \bm{\bar{x}}_{0,N}^{[N]} \, .
$$
Moreover, the Bernstein polynomial can be subdivided at $t_{\emph{div}}$ into two $N$th order Bernstein polynomials with Bernstein coefficients
$$
\bm{\bar{x}}_{0,N}^{[0]},\bm{\bar{x}}_{0,N}^{[1]},\ldots,\bm{\bar{x}}_{0,N}^{[N]} \, , \qquad \text{and} \qquad
\bm{\bar{x}}_{0,N}^{[N]},\bm{\bar{x}}_{1,N}^{[N-1]},\ldots,\bm{\bar{x}}_{N,N}^{[0]} \, .
$$
Figure \ref{app.fig:decasteljau} depicts a 2D curve defined by an $5$th order Bernstein polynomial (with Bernstein coefficients described by blue circles). The curve is subdivided into two $5$th order Bernstein polynomials, each with Bernstein coefficients described by black and red circles.

\qedprop
\end{prop}

\begin{prop}[Minimum distance] \label{app.prop:mindist}
	The minimum distance between two Bernstein polynomials $\bm{f}_N(t)$ and $\bm{g}_N(t)$, with $t\in [0,1]$, namely
	\begin{equation} \label{app.eq:mindist}
	\min_{t_a,t_b\in [0,1]} ||\bm{f}_N(t_a)-\bm{g}_N(t_b)|| \, , \qquad
	\argmin_{t_a,t_b\in [0,1]} ||\bm{f}_N(t_a)-\bm{g}_N(t_b)|| \,
	\end{equation}
	can be efficiently computed by exploiting the \emph{convex-hull} property and the de Casteljau algorithm \cite{farouki2012bernstein}, in combination with the Gilbert-Johnson-Keerthi (GJK) distance algorithm \cite{gilbert1988fast}. The latter is widely used in computer graphics and video games to compute the minimum distance between convex shapes. In \cite{chang2011computation} the authors propose an iterative procedure that uses the above tools to compute \eqref{app.eq:mindist} within a desired tolerance. This procedure is extremely useful for motion planning applications to efficiently compute the spatial clearance between two paths, or between a path and an obstacle. For example, the minimum distance between the 2D Bernstein polynomial and the point depicted in Figure \ref{app.fig:mindistpt} is computed in less than $5$ ms using an implementation in MATLAB, while the minimum distance between the 3D Bernstein polynomials depicted in Figure \ref{app.fig:mindist} is computed in less than 30 ms. The same procedure can also be employed to compute the extrema (maximum and minimum) of a Bernstein polynomial \cite{choe2017distributed}.
	
\qedprop
\end{prop}

\section{Problem formulation} \label{sec:ProblemFormulation}
This paper considers the following optimal control problem:
\begin{prob}[Problem $P$]
Determine $\bm{x}: [0,1] \to \IR^{n_x}$ and $\bm{u}: [0,1] \to \IR^{n_u}$ that minimize
\begin{equation} \label{eq:costfunc}
I(\bm{x}(t),\bm{u}(t)) = E(\bm{x}(0),\bm{x}(1))+ \int_0^{1} F(\bm{x}(t),\bm{u}(t))dt \,,
\end{equation}
subject to
\begin{align}
& \dot{\bm{x}} = \bm{f}(\bm{x}(t),\bm{u}(t))\, , \quad \forall t \in[0,1], \label{eq:dynamicconstraint} \\
& \bm{e}(\bm{x}(0),\bm{x}(1)) = \bm{0} \, , \label{eq:equalityconstraint} \\
& \bm{h}(\bm{x}(t),\bm{u}(t)) \leq \bm{0} \, , \quad \forall t\in [0,1] \label{eq:inequalityconstraint} \, ,
\end{align}
where $E:\IR^{n_x}\times \IR^{n_x}\to \IR$ and $F:\IR^{n_x}\times \IR^{n_u}\to \IR$ are the terminal and running costs, respectively, $\bm{f}:\IR^{n_x}\times \IR^{n_u}\to \IR^{n_x}$ describes the system dynamics, $\bm{e}:\IR^{n_x}\times \IR^{n_x}\to \IR^{n_e}$ is the vector of boundary conditions, and $\bm{h}:\IR^{n_x} \times \IR^{n_u}\to \IR^{n_h}$ is the vector of state and input constraints.

\qedprob
\end{prob}

\vspace{5px}
The following assumptions hold:
\begin{asm} \label{asm:Flip}
$E$, $F$, $\bm{f}$, $\bm{e}$, and $\bm{h}$ are Lipschitz continuous with respect to their arguments.
\end{asm}
\begin{asm} \label{asm:pexists}
Problem $P$ admits optimal solutions $\bm{x}^*(t)$ and $\bm{u}^*(t)$ that satisfy $\bm{x}^*(t) \in \mathcal{C}^1_{n_x}$ and $\bm{u}^*(t) \in \mathcal{C}^0_{n_u}$.
\end{asm}



\section{Bernstein approximation of Problem P} \label{sec:BernApprox}
The purpose of this section is to formulate a discretized version of Problem $P$, here referred to as Problem $P_N$, where $N$ denotes the \emph{order of approximation}.
This requires that we approximate the input and state functions, the cost function, the system dynamics, and the equality and inequality constraints in Problem $P$.

First, consider the following $N$th order vectors of Bernstein polynomials:
\begin{equation} \label{eq:polappr}
\bm{x}_N(t) =   \sum_{j=0}^N \bx_{j,N} b_{j,N}(t) \, , \quad  \bm{u}_N(t) =   \sum_{j=0}^N \bu_{j,N} b_{j,N}(t) ,
\end{equation}
with $\bm{x}_N: [0,1] \to \IR^{n_x}$, $\bm{u}_N: [0,1] \to \IR^{n_u}$, $\bx_{j,N} \in \IR^{n_x}$ and $\bu_{j,N} \in \IR^{n_u}$. Let $\bx_N \in \IR^{n_x \times (N+1)}$ and $\bu_N \in \IR^{n_u \times (N+1)}$ be defined as
$$\bx_N = [\bx_{0,N} \, , \, \ldots \, , \, \bx_{N,N}] , \qquad \bu_N = [\bu_{0,N} \, , \, \ldots \, , \, \bu_{N,N}].$$
Let $0= t_0 <t_1<\ldots<t_N=1$ be a set of equidistant \emph{time nodes}, i.e. $t_j = \frac{j}{N}$.
Then Problem $P_N$ can be stated as follows:
\begin{prob}[Problem $P_N$]
Determine $\bx_N$ and $\bu_N$ that minimize
\begin{equation} \label{eq:discCVcostfunc}
\begin{split}
& I_N(\bx_N,\bu_N) =
E(\bm{x}_N(0),\bm{x}_N(t_N))+ w \sum_{j=0}^N F(\bm{x}_N(t_j),\bm{u}_N(t_j)) \, ,
\end{split}
\end{equation}
subject to
\begin{align}
& \left\Vert \bm{\dot{x}}_N(t_j) - \bm{f}(\bm{x}_N(t_j),\bm{u}_N(t_j))\right\Vert \leq \delta_P^N \, ,
\qquad
\forall j=0,\ldots,N \,,
\label{eq:discdynamicconstraint} \\
& \bm{e}(\bm{x}_N(0),\bm{x}_N(t_N)) = \bm{0} \, , \label{eq:discequalityconstraint} \\
& \bm{h}(\bm{x}_N(t_j),\bm{u}_N(t_j)) \leq   \delta_P^N \bm{1} \, , \qquad  \forall j=0,\ldots,N  \label{eq:discinequalityconstraint} \, ,
\end{align}
where $w=\frac{1}{N+1}$, and $\delta_P^N$ is a small positive number that depends on $N$ and  converges uniformly to $0$, i.e. $\lim_{N \to \infty} \delta_P^N = 0$.

\qedprob
\end{prob}
\vspace{5px}

\begin{rem}
Compared to the constraints of Problem $P$, the dynamic and inequality constraints given by Equations \eqref{eq:discdynamicconstraint} and \eqref{eq:discinequalityconstraint} are \emph{relaxed}. Motivated by previous work on consistency of approximation theory \cite{elijah1997optimization}, the bound $\delta_P^N$, referred to as \emph{relaxation bound}, is introduced to guarantee that Problem $P_N$ has a feasible solution. As it will become clear later, the relaxation bound can be made arbitrarily small by choosing a sufficiently large order of approximation $N$.
Furthermore, note that when $N \to \infty$, then the right hand sides of Equations \eqref{eq:discdynamicconstraint} and \eqref{eq:discinequalityconstraint} equal to zero, i.e. the difference between the constraints imposed by Problems $P$ and $P_N$ vanishes.
\end{rem}

\section{Feasibility and consistency of Problem $P_N$} \label{sec:analysis}
The outcome of Problem $P_N$ is a set of optimal Bernstein coefficients $\bx^*_N$ and $\bu^*_N$ that determine the vectors of Bernstein polynomials $\bm{x}_N^*(t)$ and $\bm{u}_N^*(t)$, i.e.
\begin{equation} \label{eq:approxoptimal}
\bm{x}_N^* (t) = \sum_{j=0}^N \bx_{j,N}^* b_{j,N}(t) \, , \quad \bm{u}_N^* (t) = \sum_{j=0}^N \bu_{j,N}^* b_{j,N}(t)\, .
\end{equation}
Now we address the following theoretical issues:
\begin{enumerate}
\item existence of a feasible solution to Problem $P_{N}$,
\item convergence of the pair $(\bm{x}_N^* (t),\bm{u}_N^* (t))$ to the optimal solution of Problem $P$, given by $(\bm{x}^* (t),\bm{u}^* (t))$.
\end{enumerate}
The main results of this section are summarized in Theorems \ref{thm:existence} and \ref{thm:consistency} below.
\vspace{5px}

\begin{thm}[Feasibility] \label{thm:existence}
Let
\begin{equation} \label{eq:defdeltaPN}
\delta_P^N = C_P \max \{ W_{x^\prime}(N^{-\frac{1}{2}}) \, , \, W_{x}(N^{-\frac{1}{2}}) \, , \, W_{u}(N^{-\frac{1}{2}}) \} \, ,
\end{equation}
where $C_P$ is a positive constant independent of $N$, and $W_{x^\prime}(\cdot)$, $W_{x}(\cdot)$ and $W_{u}(\cdot)$ are the moduli of continuity of $\dot{\bm{x}}(t)$, $\bm{x}(t)$ and $\bm{u}(t)$, respectively.
Then Problem $P_N$ is feasible for arbitrary order of approximation $N \in \mathbb{Z}^+$.

\qedlem
\end{thm}

\textbf{Proof:}
Let $\bm{x}(t)$ and $\bm{u}(t)$ be a feasible solution for Problem $P$, which exists by Assumption \ref{asm:pexists}.
Let us define the Bernstein coefficients
\begin{equation} \label{eq:coeffdef}
\bx_{k,N} = \bm{x}(t_k)\, , \quad \bu_{k,N} = \bm{u}(t_k) \, , \quad \forall k \in \{0,\ldots,N\} \, ,
\end{equation}
and the resulting vectors of Bernstein polynomials as
\begin{equation} \label{eq:bezcurvesdef}
\bm{x}_N(t) = \sum_{j=0}^N \bx_{j,N} b_{j,N}(t) \, , \quad \bm{u}_N(t) = \sum_{j=0}^N \bu_{j,N} b_{j,N}(t) \, .
\end{equation}
In what follows, we show that the above polynomials satisfy the constraints in \eqref{eq:discdynamicconstraint}, \eqref{eq:discequalityconstraint} and \eqref{eq:discinequalityconstraint}, with $\delta_P^N$ defined in Equation \eqref{eq:defdeltaPN}, thus proving Theorem~\ref{thm:existence}.

Equations \eqref{eq:coeffdef} and \eqref{eq:bezcurvesdef}, together with Lemma \ref{lem:bernsteinapprox} and Assumption \ref{asm:pexists}, imply that $\bm{x}_N(t)$, $\dot{\bm{x}}_N(t)$ and $\bm{u}_N(t)$ converge uniformly to $\bm{x}(t)$, $\dot{\bm{x}}(t)$ and $\bm{u}(t)$, respectively. More precisely, for all $t \in [0,1]$ from Lemma \ref{lem:bernsteinapprox} we have
\begin{equation} \label{eq:unifconvxnprime}
\begin{split}
& ||\bm{x}_N(t) - \bm{x}(t)|| \leq C_x W_x (N^{-\frac{1}{2}}) \, ,  \\
& ||\bm{u}_N(t) - \bm{u}(t)|| \leq C_u W_u (N^{-\frac{1}{2}}) \, , \\
& ||\bm{\dot{x}}_N(t) - \bm{\dot{x}}(t)|| \leq C_{x^\prime}W_{x^{\prime}} (N^{-\frac{1}{2}}) \, ,
\end{split}
\end{equation}
 where $C_x< 5 n_x / 4$, $C_u < 5 n_u / 4$ and $C_{x^\prime} < 9n_x/4$ (see Lemma \ref{lem:bernsteinapprox}).
To prove that the dynamic constraint is satisfied, we add and subtract the term $\bm{\dot{x}}(t_k)-\bm{f}(\bm{x}(t_k),\bm{u}(t_k))$
from the left hand side of Equation \eqref{eq:discdynamicconstraint}, which yields
\begin{equation*}
\begin{split}
 || \bm{\dot{x}}_N(t_k) - \bm{f}(\bm{x}_N(t_k),\bm{u}_N(t_k))|| & \leq || \bm{\dot{x}}_N(t_k)- \bm{\dot{x}}(t_k) ||   + || \bm{\dot{x}}(t_k) - \bm{f}(\bm{x}(t_k),\bm{u}(t_k))||
\\ & \quad + ||\bm{f}(\bm{x}_N(t_k),\bm{u}_N(t_k)) - \bm{f}(\bm{x}(t_k),\bm{u}(t_k))|| \, .
\end{split}
\end{equation*}
The second term on the right hand side of the inequality above is zero (see Equation \eqref{eq:dynamicconstraint}). Using Equation \eqref{eq:unifconvxnprime} and the fact that $\bm{f}$ is Lipschitz (see Assumption \ref{asm:Flip}) with Lipschitz constant $L_f$, we get
\begin{equation*}
\begin{split}
||\bm{\dot{x}}_N(t_k) - \bm{f}(\bm{x}_N(t_k),\bm{u}_N(t_k))|| & \leq
C_{x^\prime} W_{x^{\prime}}(N^{-\frac{1}{2}})+L_f\left(C_xW_x(N^{-\frac{1}{2}})+C_u W_u(N^{-\frac{1}{2}})\right) \\
& \leq (C_{x^\prime}+L_f(C_x+C_u)) \max ( W_{x^\prime}(N^{-\delta}) , W_{x}(N^{-\delta}) , W_{u}(N^{-\delta})) \, .
\end{split}
\end{equation*}
Thus, the dynamic constraint in Equation \eqref{eq:discdynamicconstraint} is satisfied with
$\delta_P^N$ given by Equation \eqref{eq:defdeltaPN}.

Using a similar argument, the satisfaction of the constraint in Equation \eqref{eq:discinequalityconstraint} follows easily by Assumption \ref{asm:Flip}, namely that $\bm{h}$ is Lipschitz, i.e.
\begin{equation*}
\begin{split}
\bm{h}(\bm{x}_N(t_j),\bm{u}_N(t_j)) & \leq \bm{h}(\bm{x}(t_j),\bm{u}(t_j)) + ||\bm{h}(\bm{x}_N(t_j),\bm{u}_N(t_j))-\bm{h}(\bm{x}(t_j),\bm{u}(t_j))|| \\
& \leq  ||\bm{h}(\bm{x}_N(t_j),\bm{u}_N(t_j))-\bm{h}(\bm{x}(t_j),\bm{u}(t_j))|| \\
& \leq L_h (C_x + C_u) \max (W_{x}(N^{-\delta}) , W_{u}(N^{-\delta})) \, .
\end{split}
\end{equation*}
Finally, using the end point value property of Bernstein polynomials, i.e. Property \ref{app.prop:endpoint} in Section \ref{sec:MathematicalPreliminaries}, we have $\bm{x}_N(0) = \bx_{0,N}$ and $\bm{x}_N(1) = \bx_{N,N}$, which by definition implies that $\bm{e}(\bm{x}_{N}(0),\bm{x}_{N}(t_N))=\bm{e}(\bm{x}(0),\bm{x}(1))=\bm{0}$, thus proving Equation \eqref{eq:discequalityconstraint}.

\qedproof

\begin{cor} \label{cor:c2}
	If the optimal state $\bm{x}^*(t)$ and control $\bm{u}^*(t)$ solutions to Problem $P$ exist and satisfy $\dot{\bm{x}}^*(t) \in \mathcal{C}^2_{n_x}$ and $\bm{u}^*(t) \in \mathcal{C}^2_{n_u}$ in $[0,1]$, then Theorem \ref{thm:existence} holds with
	$
	\delta_{P}^N = C_PN^{-1} \, ,
	$
	where $C_P$ is a positive constant independent of $N$.
	
\qedlem
\end{cor}

\textbf{Proof:} The proof of Corollary \ref{cor:c2} follows easily by applying Lemma \ref{lem:bernsteinapproxC2} to the proof of Theorem \ref{thm:existence}.

\qedproof

\begin{rem} \label{rem:deltaPNfeas}
From the definition of $\delta_P^N$ in Theorem \ref{thm:existence} (and Corollary \ref{cor:c2}), it follows that for any arbitrarily small scalar $\epsilon_P>0$ there exists $N_1$ such that for all $N \geq N_1$, we have $\delta_P^N \leq \epsilon_P$. In other words, the relaxation bound in Problem $P_N$ can be made arbitrarily small by choosing sufficiently large $N$, while retaining the feasibility result (Theorem \ref{thm:existence} and Corollary \ref{cor:c2}).
\end{rem}

\begin{thm}[Consistency] \label{thm:consistency}
Let $\{(\bx^*_N,\bu^*_N)\}_{N=N_1}^\infty$ be a sequence of optimal solutions to Problem $P_N$, and $\{(\bm{x}_N^*(t),\bm{u}_N^*(t))\}_{N=N_1}^\infty$ be a sequence of Bernstein polynomials, given by \eqref{eq:approxoptimal}. Assume $\{(\bm{x}_N^*(t),\bm{u}_N^*(t))\}_{N=N_1}^\infty$ has a uniform accumulation point, i.e.
\begin{equation} \label{eq:convergenceasm}
\lim_{N \to \infty} (\bm{x}_N^*(t),\bm{u}_N^*(t)) = (\bm{x}^{\infty}(t),\bm{u}^{\infty}(t)) \, , \qquad \forall t \in [0,1],
\end{equation}
and assume that $\dot{\bm{x}}^\infty (t)$ and $\bm{u}^\infty (t)$ are continuous on $[0,1]$.
Then $(\bm{x}^{\infty}(t),\bm{u}^{\infty}(t))$ is an optimal solution for Problem $P$.

\qedlem
\end{thm}

\textbf{Proof:}
This proof is divided into three steps: $(1)$ we prove that $(\bm{x}^\infty(t),\bm{u}^\infty(t))$ is a feasible solution to Problem $P$; $(2)$ we show that
\begin{equation} \label{eq:proofstep2}
\lim_{N \to \infty} I_N (\bx^*_N,\bu^*_N) = I(\bm{x}^\infty(t),\bm{u}^\infty(t)) \, ;
\end{equation}
$(3)$ we prove that $(\bm{x}^\infty(t),\bm{u}^\infty(t))$ is an optimal solution of Problem $P$, i.e.
$$
I(\bm{x}^\infty(t),\bm{u}^\infty(t)) = I(\bm{x}^*(t),\bm{u}^*(t)) \, .
$$

\textbf{Step $(1)$.} First, we show that $(\bm{x}^\infty(t),\bm{u}^\infty(t))$ satisfies the dynamic constraint of Problem $P$:
$$
\bm{\dot{x}}^\infty(t)-\bm{f}(\bm{x}^\infty(t),\bm{u}^\infty(t)) = \bm{0} \, .
$$
We show this by contradiction. Assume that the above equality does not hold. Then there exists $t'$, such that
\begin{equation} \label{eq:contradiction}
||\bm{\dot{x}}^\infty(t')-\bm{f}(\bm{x}^\infty(t'),\bm{u}^\infty(t'))|| > 0 \, .
\end{equation}
Since the nodes $\{t_k\}_{k=0}^N$, $t_k = \frac{k}{N}$ are dense in $[0,1]$, there exists a sequence of indices $\{k_N\}_{N=0}^\infty$ such that
$$
\lim_{N \to \infty} t_{k_N} = t' .
$$
Then, from continuity of $\bm{\dot{x}}^\infty (t)$, $\bm{x}^\infty(t)$ and $\bm{u}^\infty(t)$, the left hand side of Equation \eqref{eq:contradiction} satisfies
\begin{equation*}
\begin{split}
& ||\bm{\dot{x}}^\infty(t')-\bm{f}(\bm{x}^\infty(t'),\bm{u}^\infty(t'))|| =
\lim_{N \to \infty} ||\bm{\dot{x}}_N^*(t_{k_N})-\bm{f}(\bm{x}_N^*(t_{k_N}),\bm{u}_N^*(t_{k_N}))|| .
\end{split}
\end{equation*}
However, the dynamic constraint in Problem $P_N$ is
$$
||\bm{\dot{x}}_N^*(t_{k_N})-\bm{f}(\bm{x}_N^*(t_{k_N}),\bm{u}_N^*(t_{k_N}))|| \leq \delta_P^N,
$$
which implies that
$$
\lim_{N \to \infty} ||\bm{\dot{x}}_N^*(t_{k_N})-\bm{f}(\bm{x}_N^*(t_{k_N}),\bm{u}_N^*(t_{k_N}))|| =  0.
$$
The above result contradicts Equation \eqref{eq:contradiction}, thus proving that $(\bm{x}^\infty(t),\bm{u}^\infty(t))$ satisfies the dynamic constraint in Equation \eqref{eq:dynamicconstraint}. The equality and inequality constraints in \eqref{eq:discequalityconstraint} and \eqref{eq:discinequalityconstraint} follow easily by an identical argument.
\vspace{5px}

\textbf{Step $(2)$.} To prove that Equation \eqref{eq:proofstep2} is satisfied we need to show the following:
\begin{equation} \label{eq:proofstep2a}
\lim_{N \to \infty} w\sum_{j=0}^{N} F(\bm{x}^*_N(t_j),\bm{u}^*_N(t_j)) =  \int_0^{1} F(\bm{x}^\infty(t),\bm{u}^\infty(t))dt \, ,
\end{equation}
and
\begin{equation} \label{eq:proofstep2b}
\lim_{N \to \infty}  E(\bm{x}^*_N(0),\bm{x}^*_N(t_N)) = E(\bm{x}^\infty(0),\bm{x}^\infty(1)) \, .
\end{equation}
Using Lemma \ref{lem:quadrature}, together with the Lipschitz assumption on $F$ (see Assumption \ref{asm:Flip}) and the continuity of $\bm{x}^\infty(t)$ and $\bm{u}^\infty(t)$, we get
$$
\lim_{N \to \infty} w \sum_{j=0}^N F(\bm{x}^\infty(t_j),\bm{u}^\infty(t_j)) = \int_0^{1} F(\bm{x}^\infty(t),\bm{u}^\infty(t))dt \, .
$$
Finally, applying the convergence assumption given by Equation \eqref{eq:convergenceasm}, the result in Equation \eqref{eq:proofstep2a} follows.
Similarly, using the Lipschitz assumption on $E$, one can show that Equation \eqref{eq:proofstep2b} holds, thus completing the proof of Step (2).
\vspace{5px}

\textbf{Step $(3)$.}
Finally, it remains to show that
$$
 I(\bm{x}^\infty(t),\bm{u}^\infty(t)) = I(\bm{x}^*(t),\bm{u}^*(t))  \, .
$$
First, recall that $\bm{x}^*(t)$ and $\bm{u}^*(t)$ are optimal solutions of Problem $P$, while $\bx^*_N$ and $\bu^*_N$ are optimal solutions of Problem $P_N$.
Let  $\tilde{\bx}_{k,N} = \bm{x}^*(t_k)$, $\tilde{\bu}_{k,N} = \bm{u}^*(t_k)$, $\forall k \in \{1,\ldots,N\}$ and
$$\tilde{\bx}_{N} = [\tilde{\bx}_{0,N} \, , \, \ldots \, , \, \tilde{\bx}_{N,N}] , \qquad \tilde{\bu}_{N} = [\tilde{\bu}_{0,N} \, , \, \ldots \, , \, \tilde{\bu}_{N,N}].$$
Following an argument similar to the one in the proof of Theorem \ref{thm:existence}, one can show that there exists $N_1$ such that for any $N\geq N_1$ the pair $(\tilde{\bx}_{N},\tilde{\bu}_{N})$ is a feasible solution of Problem $P_N$. Moreover, an argument similar to the one in the proof of Step (2) yields
\begin{equation} \label{eq:convergenceoptimalcontinuous}
\lim_{N \to \infty} I_N(\tilde{\bx}_{N},\tilde{\bu}_{N}) = I(\bm{x}^*(t),\bm{u}^*(t)) \, .
\end{equation}
Then we have
\begin{equation}
\begin{split}
& I(\bm{x}^*(t),\bm{u}^*(t)) \leq I(\bm{x}^\infty(t),\bm{u}^\infty(t))
= \lim_{N \to \infty} I_N(\bx^*_N,\bu^*_N) \leq \lim_{N \to \infty} I_N(\tilde{\bx}_{N},\tilde{\bu}_{N}) \, .
\end{split}
\end{equation}
The last inequality, combined with \eqref{eq:convergenceoptimalcontinuous}, gives
$$
I(\bm{x}^*(t),\bm{u}^*(t)) = I(\bm{x}^\infty(t),\bm{u}^\infty(t)) \, ,
$$
which completes the proof of Theorem \ref{thm:consistency}.

\qedproof

\section{Costate estimation for Problem P using Bernstein approximation} \label{sec:costate}
\subsection{First order optimality conditions of Problem P}
We start by deriving the first order necessary conditions for Problem $P$. Let ${\bm \lambda}(t): [0,1] \to \IR^{n_x}$ be the costate trajectory, and let ${\bm \mu}(t): [0,1] \to \IR^{n_h}$ and ${\bm \nu}\in \IR^{n_e}$ be the multipliers.
By defining the Lagrangian of the Hamiltonian (also known as the D-form \cite{hartl1995survey}) as
$$
\mathcal{L}(\bm{x}(t),\bm{u}(t),\bm{\lambda}(t),\bm{\mu}(t)) = \mathcal{H}(\bm{x}(t),\bm{u}(t),\bm{\lambda}(t)) + {\bm \mu}^\top (t) \bm{h}(\bm{x}(t),\bm{u}(t)) \, ,
$$
where the Hamiltonian $\mathcal{H}$ is given by 
$$
\mathcal{H}(\bm{x}(t),\bm{u}(t),\bm{\lambda}(t)) = F({\bm x}(t),{\bm u}(t)) + {\bm \lambda}^\top(t) \bm{f}(\bm{x}(t),\bm{u}(t)) \, ,
$$
the dual of Problem P can be formulated as follows  \cite{hartl1995survey}.
\begin{prob}[Problem $P_\lambda$]
Determine ${\bm x}(t)$, ${\bm u}(t)$, ${\bm \lambda}(t)$, ${\bm \mu}(t)$ and ${\bm \nu}$ that for all $t \in [0,1]$ satisfy Equations \eqref{eq:dynamicconstraint}, \eqref{eq:equalityconstraint}, \eqref{eq:inequalityconstraint} and
\begin{align}  \label{eq:Plambda1}
& \bm{\mu}^\top (t) \bm{h}({\bm x}(t),{\bm u}(t)) = 0 \, , \quad   \bm{\mu} (t) \geq 0 \, , \\ \label{eq:Plambda2}
& \dot{\bm \lambda}^\top(t) + \mathcal{L}_x (\bm{x}(t),\bm{u}(t),\bm{\lambda}(t),\bm{\mu}(t)) = \dot{\bm \lambda}^\top(t) + F_{x}({\bm x}(t),{\bm u}(t)) + {\bm \lambda}^\top(t) {\bm f}_{x} ({\bm x}(t),{\bm u}(t))
+{\bm \mu}^\top(t) {\bm h_{x}} ({\bm x}(t),{\bm u}(t)) = 0 \, , \\  \label{eq:Plambda3}
& \bm{\lambda}^\top(0) = -\bm{\nu}^\top \bm{e}_{x(0)}({\bm x}(0),{\bm x}(1)) - E_{x(0)}({\bm x}(0),{\bm x}(1)) \, , \\ \label{eq:Plambda4}
& \bm{\lambda}^\top(1) = \bm{\nu}^\top \bm{e}_{x(1)}({\bm x}(0),{\bm x}(1)) + E_{x(1)}({\bm x}(0),{\bm x}(1)) \, , \\ \label{eq:Plambda5}
& \mathcal{L}_u(\bm{x}(t),\bm{u}(t),\bm{\lambda}(t),\bm{\mu}(t)) = \bm{\lambda}^\top(t) \bm{f}_{u}({\bm x}(t),{\bm u}(t)) + F_{u}({\bm x}(t),{\bm u}(t)) + \bm{\mu}^\top(t) \bm{h}_{u}({\bm x}(t),{\bm u}(t)) = 0 \, .
\end{align}

\qedprob
\end{prob}
In the above problem, subscripts are used to denote partial derivatives, e.g. $F_{x}({\bm x},{\bm u}) = \frac{\partial}{\partial {\bm x}}F({\bm x},{\bm u})$.

The following assumptions are imposed onto Problem $P_\lambda$.
\begin{asm} \label{asm:functions}
$E,F,\bm{f},\bm{e}$ and $\bm{h}$ are continuously differentiable with respect to their arguments, and their gradients are Lipschitz continuous over the domain.
\end{asm}

\begin{asm} \label{asm:plambdaexists}
Solutions ${\bm x}^*(t)$, ${\bm u}^*(t)$, ${\bm \lambda}^*(t)$, ${\bm \mu}^*(t)$ and ${\bm \nu}^*$ of Problem $P_\lambda$ exist and satisfy ${\bm x}^*(t) \in \mathcal{C}^1_{n_x}$, $\bm{u}^*(t) \in \mathcal{C}^0_{n_u}$, ${\bm \lambda}^*(t) \in \mathcal{C}^1_{n_x}$ and $\bm{\mu}^*(t) \in \mathcal{C}^0_{n_h}$ in $[0,1]$.
\end{asm}

\begin{rem}
Notice that Problem $P_\lambda$ implicitly assumes the absence of pure state constraints in Problem $P$. If the inequality constraint in Equation \eqref{eq:inequalityconstraint} is independent of $\bm{u}(t)$, then the costate ${\bm \lambda}(t)$ must also satisfy the following jump condition \cite{hartl1995survey}:
$${\bm \lambda}(t_e^-) = {\bm \lambda}(t_e^+) + {\bm h}_{x(t_e)}^\top {\bm \eta} \, ,$$
where $t_e$ is the entry or exit time into a constrained arc in which the inequality constraint is active, $t_e^-$ and $t_e^+$ denote the left-hand side and right-hand side limits of the trajectory, respectively, and ${\bm \eta}$ is a constant covector. For simplicity, the theoretical results that will be presented in Section \ref{sec:costateresults} do not consider the jump conditions above, i.e., the inequality constraints are dependent on $\bm{u}(t)$. Nevertheless, numerical examples will be presented in Section \ref{sec:numericalresults} showing the applicability of the discretization method to pure state-constrained problems.
\end{rem}

\subsection{KKT conditions of Problem $P_N$}
Let us introduce the following $N$th order Bernstein polynomials:
\begin{equation} \label{eq:polappr_costates}
{\bm \lambda}_N(t) =   \sum_{j=0}^N \blam_{j,N} b_{j,N}(t) \, , \quad  {\bm \mu}_N(t) =   \sum_{j=0}^N \bmu_{j,N} b_{j,N}(t) ,
\end{equation}
with ${\bm \lambda}_N: [0,1] \to \IR^{n_x}$, ${\bm \mu}_N: [0,1] \to \IR^{n_h}$, $\blam_{j,N} \in \IR^{n_x}$ and $\bmu_{j,N} \in \IR^{n_h}$, and the vector $\bar{\bm{\nu}} \in \IR^{n_e}$. Finally, let
$\blam_N \in \IR^{n_x \times (N+1)}$ and $\bmu_N \in \IR^{n_u \times (N+1)}$ be defined as
$$\blam_N = [\blam_{0,N} \, , \, \ldots \, , \, \blam_{N,N}] , \qquad \bmu_N = [\bmu_{0,N} \, , \, \ldots \, , \, \bmu_{N,N}].$$
With the above notation, the Lagrangian for problem $P_N$ can be written as
\begin{equation*}
\begin{split}
\mathcal{L}_N & = E(\bm{x}_N(0),\bm{x}_N(t_N))+w \sum_{j=0}^N F({\bm x}_N(t_j),{\bm u}_N(t_j)) + \sum_{j=0}^N \bm{\lambda}_N^\top(t_j) (-{\dot{\bm x}}_N(t_j) + {\bm f}({\bm x}_N(t_j),{\bm u}_N(t_j))) \\
& \quad + \sum_{j=0}^N \bm{\mu}_N^\top(t_j) {\bm h}({\bm x}_N(t_j),{\bm u}_N(t_j)) + \bar{\bm{\nu}}^\top {\bm e}({\bm x}_N(0),{\bm x}_N(t_N))\, .
\end{split}
\end{equation*}
Then the dual of Problem $P_N$ can be stated as follows:
\begin{prob}[Problem $P_{N\lambda}$]
Determine $\bx_N$, $\bu_N$, $\blam_N$, $\bmu_N$ and $\bar{\bm \nu}$ that satisfy the primal feasibility conditions, namely Equations \eqref{eq:discdynamicconstraint}, \eqref{eq:discequalityconstraint} and \eqref{eq:discinequalityconstraint}, the complementary slackness and dual feasibility conditions
\begin{equation} \label{eq:complementarity}
\begin{split}
& \left\Vert \bm{\mu}_N^\top (t_k) \bm{h}(\bm{x}_N (t_k),\bm{u}_N (t_k)) \right\Vert  \leq N^{-1} \delta_D^N \, , \\
& \bm{\mu}_N (t_k)  \geq  - N^{-1} \delta_D^N  \bm{1} \, ,  \qquad  \forall k=0,\ldots,N \, ,
\end{split}
\end{equation}
and the stationarity conditions
\begin{equation} \label{eq:stationarity}
\left\Vert \frac{\partial \mathcal{L}_N}{\partial \bx_{k,N}} \right\Vert \leq \delta_D^N \, , \quad \left\Vert \frac{\partial \mathcal{L}_N}{\partial \bu_{k,N}} \right\Vert \leq \delta_D^N \, , \quad \forall k = 0,\ldots,N ,
\end{equation}
where $\delta_D^N$ is a small positive number that depends on $N$ and  satisfies $\lim_{N \to \infty} \delta_D^N = 0$.

\qedprob
\end{prob}

At this point one might expect results similar to the ones in Section \ref{sec:analysis}, i.e. feasibility (Theorem \ref{thm:existence}) and consistency (Theorem \ref{thm:consistency}). Nevertheless, similarly to most results on costate estimation  \cite{garg2011direct,gong2008connections,singh2009verification}, this is not  the case, and additional conditions must be added to Equations \eqref{eq:discdynamicconstraint}-\eqref{eq:discinequalityconstraint}, \eqref{eq:complementarity} and \eqref{eq:stationarity} in order to obtain consistent approximations of the solutions of Problem $P_{\lambda}$. These conditions, often referred to as \emph{closure conditions} in the literature, are given as follows:
\begin{align}
& \label{eq:closurecond1} \left\Vert \frac{\bm{\lambda}_N^\top(0)}{w}+\bar{\bm \nu}^\top {\bm e}_{x(0)}(\bm{x}_N(0),\bm{x}_N(t_N)) + E_{x(0)}(\bm{x}_N(0),\bm{x}_N(t_N))\right\Vert \leq \delta_D^N  \\
& \label{eq:closurecond2} \left\Vert \frac{\bm{\lambda}_N^\top(t_N)}{w} - \bar{\bm \nu}^\top {\bm e}_{x(1)}(\bm{x}_N(0),\bm{x}_N(t_N))- E_{x(1)}(\bm{x}_N(0),\bm{x}_N(t_N)) \right\Vert  \leq \delta_D^N  \, .
\end{align}
In other words, the closure conditions are constraints that must be added to Problem $P_{N\lambda}$ so that the solution of this problem approximates the solution of Problem $P_\lambda$. We notice that the conditions given above are discrete approximations of the conditions given by Equations \eqref{eq:Plambda3} and \eqref{eq:Plambda4}.
With this setup, we define the following problem:
\begin{prob}[Problem $P_{N\lambda}^{clos}$]
Determine $\bx_N$, $\bu_N$, $\blam_N$, $\bmu_N$ and $\bar{\bm \nu}$ that satisfy the primal feasibility conditions, namely Equations \eqref{eq:discdynamicconstraint}, \eqref{eq:discequalityconstraint} and \eqref{eq:discinequalityconstraint}, the complementary slackness and dual feasibility conditions \eqref{eq:complementarity}, the stationarity conditions \eqref{eq:stationarity}, and the closure conditions \eqref{eq:closurecond1} and \eqref{eq:closurecond2}.

\qedprob
\end{prob}

The solution of Problem $P_{N\lambda}^{clos}$ presents a set of optimal Bernstein coefficients $\bx^*_N$, $\bu^*_N$, $\blam^*_N$, $\bmu^*_N$ (which determine the Bernstein polynomials $\bm{x}_N^*(t)$, $\bm{u}_N^*(t)$, $\bm{\lambda}_N^*(t)$ and $\bm{\mu}_N^*(t)$) and a vector $\bar{\bm \nu}^*$.

\section{Feasibility and consistency of Problem $P_{N \lambda}^{clos}$} \label{sec:costateresults}
The objective of this section is to investigate the ability of the solutions of Problem $P_{N\lambda}^{clos}$ to approximate the solutions of Problem $P_\lambda$.
In what follows, we first show the existence of a solution for Problem $P_{N\lambda}^{clos}$ that satisfies also the closure conditions (feasibility).
Second, we investigate the convergence properties of this solution as $N \to \infty$ (consistency). Third, by combining these two results, we finally formulate the \emph{covector mapping theorem} for Bernstein approximations, which provides a bijective map (covector mapping) between the solution of Problem $P_{N\lambda}^{clos}$ and the solution of Problem $P_\lambda$.
The main results of this section are summarized in the three theorems below.

\begin{thm}[Feasibility] \label{thm:existence_dual}
Let
\begin{equation} \label{eq:defdeltaDN}
\delta_D^N = C_D \max \{ \delta_P^N \, , \, W_{\lambda^\prime}(N^{-\frac{1}{2}})\, , \, W_{\lambda}(N^{-\frac{1}{2}}) \, , \, W_{\mu}(N^{-\frac{1}{2}}) \} \, ,
\end{equation}
where $C_D$ is a positive constant independent of $N$, $\delta_P^N $ was defined in Equation \eqref{eq:defdeltaPN}, and $W_{\lambda^\prime}(\cdot)$, $W_{\lambda}(\cdot)$, and $W_{\mu}(\cdot)$ are the moduli of continuity of $\dot{\bm{\lambda}}(t)$, $\bm{\lambda}(t)$ and $\bm{\mu}(t)$, respectively.
Then Problem $P_{N\lambda}^{clos}$ is feasible for arbitrary order of approximation $N \in \mathbb{Z}^+$.

\qedlem
\end{thm}

\textbf{Proof:} Similar to the proof of Theorem \ref{thm:existence}, this proof follows by constructing a solution for Problem $P_{N\lambda}^{clos}$, with $\delta_D^N$ given by Equation \eqref{eq:defdeltaDN}.
To this end, let $x(t)$, $u(t)$, $\lambda(t)$, $\mu(t)$ and $\nu$ be a solution of Problem $P_\lambda$, which exists by Assumption \ref{asm:plambdaexists}, and define
\begin{equation} \label{eq:defxu}
\bx_{j,N} = \bm{x}(t_j) \, , \quad \bu_{j,N} = \bm{u}(t_j) \, ,
\end{equation}
\begin{equation} \label{eq:deflammu}
{\bm{\bar{\lambda}}}_{j,N} = w \bm{\lambda}(t_j) \, , \quad  {\bm{\bar{\mu}}}_{j,N} = w \bm{\mu}(t_j) \, , \quad \bm{\bar{\nu}} = \bm{\nu} \, ,
\end{equation}
for all $j = 0 , \ldots , N$, $t_j = \frac{j}{N}$, $w = \frac{1}{N+1}$, with corresponding Bernstein polynomials given by
\begin{equation} \label{eq:bezcurvesdef_dual}
\bm{x}_N(t) = \sum_{j=0}^N \bx_{j,N} b_{j,N}(t) \, , \quad \bm{u}_N(t) = \sum_{j=0}^N \bu_{j,N} b_{j,N}(t) \, , \quad {\bm \lambda}_N(t) = \sum_{j=0}^N {{\blam}}_{j,N} b_{j,N}(t) \, , \quad  {\bm \mu}_N(t) =   \sum_{j=0}^N {{\bmu}}_{j,N} b_{j,N}(t) .
\end{equation}
The remainder of this proof shows that $\bm{x}_{N}(t)$, $\bm{u}_{N}(t)$, ${{\bm{\lambda}}}_{N}(t)$ , ${{\bm{\mu}}}_{N}(t)$ and $\bar{\bm \nu}$ given above satisfy Equations \eqref{eq:complementarity}-\eqref{eq:closurecond2} (we notice that the satisfaction of Equations \eqref{eq:discdynamicconstraint}-\eqref{eq:discinequalityconstraint} has already been addressed in the proof of Theorem \ref{thm:existence}). We start by defining the Bernstein coefficients $\tilde{\bar{\bm \lambda}}_{j,N}$ and $\tilde{\bar{\bm \mu}}_{j,N}$ as follows
\begin{equation} \label{eq:deflammutilde}
\tilde{\bar{\bm \lambda}}_{j,N} = \frac{{\bar{\bm \lambda}}_{j,N}}{w} \, , \quad \tilde{\bar{\bm \mu}}_{j,N} = \frac{{\bar{\bm \mu}}_{j,N}}{w} \, ,
\end{equation}
with corresponding Bernstein polynomials given by
\begin{equation*}
\tilde{\bm \lambda}_N(t) =   \sum_{j=0}^N \tilde{\bar{\bm{\lambda}}}_{j,N} b_{j,N}(t) \, , \quad  \tilde{\bm \mu}_N(t) =   \sum_{j=0}^N \tilde{\bar{\bm{\mu}}}_{j,N} b_{j,N}(t) .
\end{equation*}
Notice that
\begin{equation} \label{eq:transformationtilde}
\tilde{\bm \lambda}_N(t) = \frac{{\bm \lambda}_N(t)}{w} \, , \qquad \tilde{\bm \mu}_N(t) = \frac{{\bm \mu}_N(t)}{w} \, .
\end{equation}
Combining Equations \eqref{eq:defxu}, \eqref{eq:deflammu} and \eqref{eq:deflammutilde} and using Assumption \ref{asm:plambdaexists} and Lemma \ref{lem:bernsteinapprox}, we get
\begin{equation} \label{eq:boundedcostates}
\begin{split}
& ||{\bm{x}}_N(t)-\bm{x}(t)|| \leq C_x W_{x}(N^{-\frac{1}{2}}) \, , \quad  ||{\bm{u}}_N(t)-\bm{u}(t)|| \leq C_u W_{u}(N^{-\frac{1}{2}}) \, , \quad ||\dot{\bm{x}}_N(t)-\dot{\bm{x}}(t)|| \leq C_{x^\prime} W_{{x^\prime}}(N^{-\frac{1}{2}}) \, ,  \\
& ||\tilde{\bm{\lambda}}_N(t)-\bm{\lambda}(t)|| \leq C_\lambda W_{\lambda}(N^{-\frac{1}{2}}) \, , \quad  ||\tilde{\bm{\mu}}_N(t)-\bm{\mu}(t)|| \leq C_\mu W_{\mu}(N^{-\frac{1}{2}}) \, , \quad ||\dot{\tilde{\bm{\lambda}}}_N(t)-\dot{\bm{\lambda}}(t)|| \leq C_{\lambda^\prime} W_{\lambda^\prime}(N^{-\frac{1}{2}}) \, ,
\end{split}
\end{equation}
where $C_\lambda < \frac{5 n_x}{4}  \, , \quad C_\mu < \frac{5n_h}{4} \, , \ C_{\lambda^\prime} < \frac{9n_x}{4}$ and $W_{\lambda}(\cdot)$, $W_{\mu}(\cdot)$ and $W_{\mu}(\cdot)$ are the moduli of continuity of $\bm{\lambda} (t)$, $\bm{\mu} (t)$ and $\dot{\bm{\lambda}} (t)$, respectively.

Now we show that the bound in Equation \eqref{eq:complementarity} is satisfied. Using Equation \eqref{eq:transformationtilde}, and adding and subtracting $w({\bm{\mu}}^\top (t_k) \bm{h}(\bm{x}_N (t_k),\bm{u}_N (t_k)) + {\bm{\mu}}^\top (t_k) \bm{h}(\bm{x} (t_k),\bm{u} (t_k)))$, we get
\begin{equation*}
\begin{split}
\Vert \bm{\mu}_N^\top (t_k) \bm{h}(\bm{x}_N (t_k),\bm{u}_N (t_k)) \Vert & = \Vert w \tilde{\bm{\mu}}_N^\top (t_k) \bm{h}(\bm{x}_N (t_k),\bm{u}_N (t_k))|| \\
&  \leq w \Vert (\tilde{\bm{\mu}}_N^\top (t_k) - {\bm{\mu}}^\top (t_k)) \bm{h}(\bm{x}_N (t_k),\bm{u}_N (t_k)) \Vert + w \Vert {\bm{\mu}}^\top (t_k) \bm{h}(\bm{x}(t_k),\bm{u}(t_k)) \Vert \\
&  \quad + w \Vert {\bm{\mu}}^\top (t_k) (\bm{h}(\bm{x}_N (t_k),\bm{u}_N (t_k))-\bm{h}(\bm{x}(t_k),\bm{u}(t_k))\Vert
\end{split}
\end{equation*}
Using Equation \eqref{eq:Plambda1}, the above inequality reduces to
\begin{equation*}
\begin{split}
\Vert \bm{\mu}_N^\top (t_k) \bm{h}(\bm{x}_N (t_k),\bm{u}_N (t_k)) \Vert & \leq
w \Vert (\tilde{\bm{\mu}}_N^\top (t_k) - {\bm{\mu}}^\top (t_k)) \bm{h}(\bm{x}_N (t_k),\bm{u}_N (t_k)) \Vert \\
& \quad + w \Vert {\bm{\mu}}^\top (t_k) (\bm{h}(\bm{x}_N (t_k),\bm{u}_N (t_k))-\bm{h}(\bm{x}(t_k),\bm{u}(t_k))\Vert  \\
&  \leq w ||\bm{h}(\bm{x}_N (t_k),\bm{u}_N (t_k))|| C_\mu W_{\mu}(N^{-\frac{1}{2}}) \\
&  \quad + w \Vert {\bm{\mu}}^\top (t_k) \Vert L_h (C_x W_x(N^{-\frac{1}{2}})+C_u W_u(N^{-\frac{1}{2}})), \\
\end{split}
\end{equation*}
where we used the bounds in Equation \eqref{eq:boundedcostates} together with the Lipschitz assumption on $\bm{h}$ (see Assumptions \ref{asm:functions}). Finally, from using Assumptions \ref{asm:functions} and \ref{asm:plambdaexists}, it follows that $\bm{h}$ and ${\bm{\mu}}$ are bounded on $[0,1]$ with bounds $h_{\max}$ and $\mu_{\max}$, respectively. Therefore, we get
\begin{equation*}
\begin{split}
\Vert \bm{\mu}_N^\top (t_k) \bm{h}(\bm{x}_N (t_k),\bm{u}_N (t_k)) \Vert \leq  w [h_{\max} C_\mu W_{\mu}(N^{-\frac{1}{2}})+\mu_{\max} L_h (C_x W_x(N^{-\frac{1}{2}})+C_u W_u(N^{-\frac{1}{2}}))] ,
\end{split}
\end{equation*}
which implies that the bound in Equation \eqref{eq:complementarity} is satisfied with $\delta_D^N$ given by Equation \eqref{eq:defdeltaDN} and $C_D > h_{\max} C_\mu + \mu_{\max}L_h(C_x+C_u)$.
Similarly,
\begin{equation*}
\begin{split}
& \bm{\mu}_N(t_k) = w \tilde{\bm \mu}_N (t_k) \geq w {\bm \mu}(t_k)-w\Vert {\bm \mu(t_k)} - \tilde{\bm \mu}_N(t_k) \Vert \bm{1}
\geq - N^{-1} C_{\mu} W_{\mu} (N^{-\frac{1}{2}}) \bm{1} \, ,
\end{split}
\end{equation*}
which proves that Equation \eqref{eq:complementarity} holds.

Now consider the left equation in \eqref{eq:stationarity}. For $k=0$ we have
\begin{equation} \label{eq:comp1}
\begin{split}
\left\Vert \frac{\partial \mathcal{L}_N}{\partial {\bx}_{0,N}} \right\Vert & = \left\Vert E_{x(0)} (\bm{x}_N(0),\bm{x}_N (t_N))
+ w \sum_{j = 0}^N F_x(\bm{x}_N(t_j),\bm{u}_N (t_j)) b_{0,N}(t_j) \right.  \\
& \quad + \sum_{j = 0}^N \bm{\lambda}^\top_N(t_j)\left[-\dot{b}_{0,N}(t_j)+{\bm f}_x(\bm{x}_N(t_j),\bm{u}_N (t_j)) b_{0,N}(t_j)\right] \\
& \quad \left. + \sum_{j = 0}^N \bm{\mu}^\top_N(t_j){\bm h}_x(\bm{x}_N(t_j),\bm{u}_N (t_j))b_{0,N}(t_j) + \bar{\bm \nu}^\top {\bm e}_{x(0)}(\bm{x}_N(0),\bm{x}_N (t_N)) \right\Vert  \, .
\end{split}
\end{equation}
Substituting $w\tilde{\bm \lambda}_N(t_j)={\bm \lambda}_N(t_j)$ and $w\tilde{\bm \mu}_N(t_j)={\bm \mu}_N(t_j)$, the equation above can be written as
\begin{equation} \label{eq:comp12}
\begin{split}
\left\Vert \frac{\partial \mathcal{L}_N}{\partial {\bx}_{0,N}} \right\Vert & = \left\Vert E_{x(0)} (\bm{x}_N(0),\bm{x}_N (t_N))
+ w \sum_{j = 0}^N F_x(\bm{x}_N(t_j),\bm{u}_N (t_j)) b_{0,N}(t_j) \right.  \\
& \quad + w \sum_{j = 0}^N \tilde{\bm{\lambda}}^\top_N(t_j)\left[-\dot{b}_{0,N}(t_j)+{\bm f}_x(\bm{x}_N(t_j),\bm{u}_N (t_j)) b_{0,N}(t_j)\right] \\
& \quad \left. + w\sum_{j = 0}^N \tilde{\bm{\mu}}^\top_N(t_j){\bm h}_x(\bm{x}_N(t_j),\bm{u}_N (t_j))b_{0,N}(t_j) + \bar{\bm \nu}^\top {\bm e}_{x(0)}(\bm{x}_N(0),\bm{x}_N (t_N)) \right\Vert  \, .
\end{split}
\end{equation}
Notice that the following inequalities are satisfied:
\begin{subequations} \label{eq:inequalities1}
\begin{eqnarray}
\label{eq:inequalities11}
\left\Vert w \sum_{j = 0}^N F_x(\bm{x}_N(t_j),\bm{u}_N (t_j)) b_{0,N}(t_j) - \int_0^1 F_x(\bm{x}(t),\bm{u}(t)) b_{0,N}(t) dt \right\Vert \leq \bar{C}_1 (N^{-\frac{1}{2}} + W_{x}(N^{-\frac{1}{2}})+W_{u}(N^{-\frac{1}{2}})) \\[10pt]
\label{eq:inequalities12}
\left\Vert  w \sum_{j = 0}^N \tilde{\bm{\lambda}}^\top_N(t_j)\dot{b}_{0,N}(t_j) - \int_0^1 {\bm{\lambda}}^\top(t)\dot{b}_{0,N}(t)dt \right\Vert \leq \bar{C}_2 (N^{-\frac{1}{2}} +  W_{\lambda}(N^{-\frac{1}{2}})) \\[10pt]
\label{eq:inequalities13}
\left\Vert  w \sum_{j = 0}^N \tilde{\bm{\lambda}}^\top_N(t_j){\bm f}_x(\bm{x}_N(t_j),\bm{u}_N (t_j)) b_{0,N}(t_j) - \int_0^1 {\bm{\lambda}}^\top(t){\bm f}_x(\bm{x}(t),\bm{u} (t)) b_{0,N}(t)dt \right\Vert \leq \qquad \qquad \qquad \\
\qquad \bar{C}_3 (N^{-\frac{1}{2}} +  W_{\lambda}(N^{-\frac{1}{2}})+W_{x}(N^{-\frac{1}{2}})+W_{u}(N^{-\frac{1}{2}})) \nonumber \\[10pt]
\label{eq:inequalities14}
\left\Vert  w\sum_{j = 0}^N \tilde{\bm{\mu}}^\top_N(t_j){\bm h}_x(\bm{x}_N(t_j),\bm{u}_N (t_j))b_{0,N}(t_j) -  \int_0^1 {\bm{\mu}}^\top(t){\bm h}_x(\bm{x}(t),\bm{u} (t))b_{0,N}(t) dt  \right\Vert \leq \qquad \qquad \qquad \\ \nonumber
\qquad \bar{C}_4 (N^{-\frac{1}{2}} +  W_{\mu}(N^{-\frac{1}{2}})+W_{x}(N^{-\frac{1}{2}})+W_{u}(N^{-\frac{1}{2}})) \, ,
\end{eqnarray}
\end{subequations}
for some positive $\bar{C}_1$, $\bar{C}_2$, $\bar{C}_3$ and $\bar{C}_4$ independent of $N$. A proof of the above inequalities is given in Appendix \ref{app:inequalities1}.
Then the combination of Equations \eqref{eq:comp12} and \eqref{eq:inequalities1} yields the following inequality
\begin{equation} \label{eq:comp13}
\begin{split}
\left\Vert \frac{\partial \mathcal{L}_N}{\partial {\bx}_{0,N}} \right\Vert & \leq \left\Vert E_{x(0)} (\bm{x}(0),\bm{x} (1))
+ \int_0^1 F_x(\bm{x}(t),\bm{u}(t)) b_{0,N}(t) dt - \int_0^1 {\bm{\lambda}}^\top(t)\dot{b}_{0,N}(t)dt \right. \\
& \quad + \int_0^1 {\bm{\lambda}}^\top(t){\bm f}_x(\bm{x}(t),\bm{u} (t)) b_{0,N}(t)dt
 \left. + \int_0^1 {\bm{\mu}}^\top(t){\bm h}_x(\bm{x}(t),\bm{u} (t))b_{0,N}(t) dt + \bar{\bm \nu}^\top {\bm e}_{x(0)}(\bm{x}(0),\bm{x} (1)) \right\Vert  \\
& \quad + \bar{C} \max \left\{ N^{-\frac{1}{2}},W_{x}(N^{-\frac{1}{2}}),W_{u}(N^{-\frac{1}{2}}),W_{\lambda}(N^{-\frac{1}{2}}),W_{\mu}(N^{-\frac{1}{2}}) \right\} \, ,
\end{split}
\end{equation}
with $\bar{C} \geq 4 \max\{\bar{C}_1,\bar{C}_2,\bar{C}_3,\bar{C}_4\}$.
Using integration by parts, we have $\int_0^{1} {\bm{\lambda}}^\top(t)\dot{b}_{0,N}(t)dt = - \int_0^{1} \dot{{\bm{\lambda}}}^\top(t)b_{0,N}(t)dt + [{\bm{\lambda}}^\top(t)b_{0,N}(t)]_{0}^{1}$. Thus, since $b_{0,N}(0)=1, b_{N,N}(0)=0$, the above inequality becomes
\begin{equation} \label{eq:comp3}
\begin{split}
\left\Vert \frac{\partial \mathcal{L}_N}{\partial {\bx}_{0,N}} \right\Vert & \leq \left\Vert E_{x(0)} (\bm{x}(0),\bm{x}(1)) + {\bm{\lambda}}^\top(0) + {\bm \nu}^\top {\bm e}_{x(0)}(\bm{x}(0),\bm{x} (1)) \right. \\
& \quad \left. + \int_0^{1} \left[ \dot{{\bm{\lambda}}}^\top(t) + F_x(\bm{x}(t),\bm{u}(t)) + {\bm{\lambda}}^\top(t){\bm f}_x(\bm{x}(t),\bm{u} (t)) + {\bm{\mu}}^\top(t){\bm h}_x(\bm{x}(t),\bm{u} (t)) \right] b_{0,N}(t)dt \right\Vert \\
& \quad + \bar{C} \max \left\{ N^{-\frac{1}{2}},W_{x}(N^{-\frac{1}{2}}),W_{u}(N^{-\frac{1}{2}}),W_{\lambda}(N^{-\frac{1}{2}}),W_{\mu}(N^{-\frac{1}{2}}) \right\}  \, .
\end{split}
\end{equation}
Finally, using Equations \eqref{eq:Plambda2} and \eqref{eq:Plambda3}, the above inequality reduces to the left condition in Equation \eqref{eq:stationarity} for $k = 0$, with $\delta_D^N$ given by Equation \eqref{eq:defdeltaDN} and $C_D \geq \bar{C}$.
The same condition for $k=1,\ldots,N$ can be shown to be satisfied using an identical argument. The stationarity condition in the right of Equation \eqref{eq:stationarity}
can also be verified similarly, and the computations are thus omitted. To show that the closure condition \eqref{eq:closurecond1} is satisfied we use the definitions in Equations \eqref{eq:defxu} and \eqref{eq:deflammu} together with the end point values property of Bernstein polynomials, Property \ref{app.prop:endpoint} in Section \ref{sec:MathematicalPreliminaries}, which gives
\begin{align*}
& \left\Vert \frac{\bm{\lambda}_N^\top(0)}{w}+\bar{\bm \nu}^\top {\bm e}_{x(0)}(\bm{x}_N(0),\bm{x}_N(t_N)) + E_{x(0)}(\bm{x}_N(0),\bm{x}_N(t_N))\right\Vert \leq  \\
& \qquad \qquad \left\Vert \bm{\lambda}^\top(0)+{\bm \nu}^\top {\bm e}_{x(0)}(\bm{x}(0),\bm{x}(1)) + E_{x(0)}(\bm{x}(0),\bm{x}(1))\right\Vert = 0  \, ,
\end{align*}
where the last equality follows from Equation \eqref{eq:Plambda3}. An identical argument can be used to show that the closure condition \eqref{eq:closurecond2} holds, thus completing the proof of Theorem \ref{thm:existence_dual}.

\qedproof

\begin{cor} \label{cor:c2_dual}
	If solutions ${\bm x}^*(t)$, ${\bm u}^*(t)$, ${\bm \lambda}^*(t)$, ${\bm \mu}^*(t)$ and ${\bm \nu}^*$ of Problem $P_\lambda$ exist and satisfy $\dot{\bm x}^*(t) \in \mathcal{C}^2_{n_x}$, $\bm{u}^*(t) \in \mathcal{C}^2_{n_u}$, $\dot{\bm \lambda}^*(t) \in \mathcal{C}^2_{n_x}$, and $\bm{\mu}^*(t) \in \mathcal{C}^2_{n_h}$ in $[0,1]$, then Theorem \ref{thm:existence_dual} holds with
	$ \delta_{P}^N = C_P N^{-1} $ and $ \delta_{D}^N = C_D N^{-1} \, ,$
	where $C_P$ and $C_D$ are positive constants independent of $N$.

\qedlem	
\end{cor}

\textbf{Proof:} The proof of Corollary \ref{cor:c2_dual} follows easily by applying Lemma \ref{lem:bernsteinapproxC2} to the proof of Theorem \ref{thm:existence_dual}.

\qedproof

\begin{rem}
Similarly to Remark \ref{rem:deltaPNfeas}, for arbitrarily small scalar $\epsilon_D>0$ there exists $N_1$ such that for all $N \geq N_1$, we have $\delta_D^N \leq \epsilon_D$; i.e., the relaxation bound in Problem $P_{N \lambda}^\text{clos}$ can be made arbitrarily small by choosing sufficiently large $N$.
\end{rem}

\begin{thm}[Consistency] \label{thm:consistency_dual}
Let $\{({\bx}_N^*,{\bu}_u^*, {\blam}_N^*,{\bmu}_N^*,\bar{\bm \nu}^*)\}_{N=N_1}^\infty$ be a sequence of solutions of Problem $P_{N\lambda}^{\text{clos}}$. Consider the sequence of transformed solutions
${\{({\bx}_N^*,{\bu}_N^*, \tilde{\blam}_N^*,\tilde{\bmu}_N^*,\bar{\bm \nu}^*)\}_{N=N_1}^\infty}$, with
$$
\tilde{\blam}_{j,N}^* = \frac{{\blam}_{j,N}^*}{w} \, , \quad \tilde{\bmu}_{j,N}^* = \frac{{\blam}_{j,N}^*}{w} \, ,
$$
and the corresponding polynomial approximation $\{({\bm x}_N^*(t),{\bm u}_N^*(t), \tilde{\bm \lambda}_N^*(t),\tilde{\bm \mu}_N^*(t),\bar{\bm \nu}^*)\}_{N=N_1}^\infty$. Assume that the latter has a uniform accumulation point, i.e.
\begin{equation*}
\begin{split}
& \lim_{N \to \infty} ({\bm x}_N^*(t),{\bm u}_N^*(t), \tilde{\bm \lambda}_N^*(t),\tilde{\bm \mu}_N^*(t),\bar{\bm \nu}^*) =  ({\bm x}^\infty(t),{\bm u}^\infty(t), \tilde{\bm \lambda}^\infty(t),\tilde{\bm \mu}^\infty(t),\bar{\bm \nu}^\infty) \, , \qquad \forall t \in [0,1] ,
\end{split}
\end{equation*}
and assume $\dot{\bm x}^\infty(t)$, ${\bm u}^\infty(t)$, $\dot{\tilde{\bm \lambda}}^\infty(t)$ and $\tilde{\bm \mu}^\infty(t)$ are continuous on $[0,1]$. Then, $({\bm x}^\infty(t),{\bm u}^\infty(t), \tilde{\bm \lambda}^\infty(t),\tilde{\bm \mu}^\infty(t),\bar{\bm \nu}^\infty)$ is a solution of Problem $P_\lambda$.

\qedlem
\end{thm}

\textbf{Proof:} The objective is to show that ${\bm x}^\infty(t),{\bm u}^\infty(t), \tilde{\bm \lambda}^\infty(t),\tilde{\bm \mu}^\infty(t)$ and $\bar{\bm \nu}^\infty$ satisfy Equations \eqref{eq:dynamicconstraint}-\eqref{eq:inequalityconstraint} and \eqref{eq:Plambda1}-\eqref{eq:Plambda5}. The satisfaction of Equations \eqref{eq:dynamicconstraint}-\eqref{eq:inequalityconstraint} has been demonstrated in the proof of Theorem \ref{thm:consistency}. We start by showing Equation \eqref{eq:Plambda1}, and we do so using a proof by contradiction. Assume that ${\bm x}^\infty(t),{\bm u}^\infty(t), \tilde{\bm \mu}^\infty(t)$ do not satisfy Equation \eqref{eq:Plambda1}. Then there exists $t' \in [0,1]$, such that
\begin{equation} \label{eq:contradict}
\Vert {\tilde{\bm{\mu}}}^{\infty \top} (t') \bm{h}({\bm x}^\infty(t'),{\bm u}^\infty(t')) \Vert > 0  .
\end{equation}
Since the nodes $\{t_k\}_{k=0}^N$ are dense in $[0,1]$, there exists a sequence of indices $\{k_N\}_{N=0}^\infty$ such that
$$
\lim_{N \to \infty} t_{k_N} = t' ,
$$
which implies
$$
\lim_{N \to \infty} \Vert \tilde{\bm \mu}^{\infty} (t')-\tilde{\bm \mu}^{\infty} (t_{k_N}) \Vert = 0 \, ,
$$
$$
\lim_{N \to \infty} \Vert {\bm x}^{\infty} (t')-{\bm x}^{\infty} (t_{k_N}) \Vert = 0 \, ,
$$
$$
\lim_{N \to \infty} \Vert {\bm u}^{\infty} (t')-{\bm u}^{\infty} (t_{k_N}) \Vert = 0 \, .
$$
Then we have
\begin{equation*}
\begin{split}
 || {\tilde{\bm \mu}}^{\infty \top}(t') \bm{h}(\bm{x}^\infty(t'),\bm{u}^\infty(t'))  || & \leq \lim_{N \to \infty} || ({\tilde{\bm \mu}}_N^{* \top}(t')-{\tilde{\bm \mu}}_N^{* \top}(t_{k_N}))\bm{h}(\bm{x}_N^*(t'),\bm{u}_N^*(t')) || \\
& \qquad  + \lim_{N \to \infty} || {\tilde{\bm \mu}}_N^{* \top}(t_{k_N})(\bm{h}(\bm{x}_N^*(t'),\bm{u}_N^*(t')) -\bm{h}(\bm{x}_N^*(t_{k_N}),\bm{u}_N^*(t_{k_N}))) || \\
& \qquad  + \lim_{N \to \infty} || {\tilde{\bm \mu}}_N^{* \top}(t_{k_N})\bm{h}(\bm{x}_N^*(t_{k_N}),\bm{u}_N^*(t_{k_N})) || \\
&  = \lim_{N \to \infty} \frac{1}{w} || {{\bm \mu}}_N^{* \top}(t_{k_N})\bm{h}(\bm{x}_N^*(t_{k_N}),\bm{u}_N^*(t_{k_N})) ||  = 0 \, ,
\end{split}
\end{equation*}
where we used Equation \eqref{eq:complementarity}. This contradicts Equation \eqref{eq:contradict}. Similarly, we can show that $\tilde{\bm \mu}^\infty(t) \geq 0$, thus proving that ${\bm x}^\infty(t),{\bm u}^\infty(t)$ and $\tilde{\bm \mu}^\infty(t)$ satisfy Equation \eqref{eq:Plambda1}.

Furthermore, we notice that if ${\bm x}^\infty(t),{\bm u}^\infty(t), \tilde{\bm \lambda}^\infty(t),\tilde{\bm \mu}^\infty(t)$ and $\bar{\bm \nu}^\infty$ satisfy Equations \eqref{eq:stationarity}-\eqref{eq:closurecond2}, then the following holds for all $k=0,\ldots,N$:
$$\left\Vert {\tilde{\bm{\lambda}}}^{\infty \top} (0)+{\bar{\bm \nu}}^{\infty \top} {\bm e}_{x(0)}(\bm{x}^\infty(0),\bm{x}^\infty(1)) + E_{x(0)}(\bm{x}^{\infty}(0),\bm{x}^\infty(1))\right\Vert = 0 \, , $$
$$\left\Vert {\bm{\lambda}}^{\infty \top}(1)-{\bar{\bm \nu}}^{\infty \top} {\bm e}_{x(1)}(\bm{x}^\infty(0),\bm{x}^\infty(1)) - E_{x(1)}(\bm{x}^\infty(0),\bm{x}^\infty(1)) \right\Vert  = 0 \nonumber \, , $$
$$\left\Vert \int_0^{1} \left[ {\dot{\tilde{\bm{\lambda}}}}^{\infty \top}(t) + F_x(\bm{x}^\infty(t),\bm{u}^\infty (t))+ {\tilde{\bm{\lambda}}}^{\infty \top}(t){\bm f}_x(\bm{x}^\infty(t),\bm{u}^\infty (t)) + {\tilde{\bm{\mu}}}^{\infty \top}(t){\bm h}_x(\bm{x}^\infty(t),\bm{u}^\infty(t)) \right] b_{k,N}(t)dt \right\Vert = 0 \, , $$
$$\left\Vert \int_0^{1} \left[ F_u(\bm{x}^\infty(t),\bm{u}^\infty (t))+ {\tilde{\bm{\lambda}}}^{\infty \top}(t){\bm f}_u(\bm{x}^\infty(t),\bm{u}^\infty (t)) + {\tilde{\bm{\mu}}}^{\infty \top}(t){\bm h}_u(\bm{x}^\infty(t),\bm{u}^\infty(t)) \right] b_{k,N}(t)dt \right\Vert = 0 \, . $$
 Since $\{b_{k,N}(t)\}_{k=0}^N$ is a linearly independent basis set, the last two equations above imply
\begin{equation*}
\begin{split}
& \left\Vert {\dot{\tilde{\bm{\lambda}}}}^{\infty \top}(t) + F_x(\bm{x}^\infty(t),\bm{u}^\infty (t)) + {\tilde{\bm{\lambda}}}^{\infty \top}(t){\bm f}_x(\bm{x}^\infty(t),\bm{u}^\infty (t)) +  {\tilde{\bm{\mu}}}^{\infty \top}(t){\bm h}_x(\bm{x}^\infty(t),\bm{u}^\infty (t)) \right\Vert = 0 \, ,
\end{split}
\end{equation*}
\begin{equation*}
\begin{split}
& \left\Vert F_u(\bm{x}^\infty(t),\bm{u}^\infty (t)) + {\tilde{\bm{\lambda}}}^{\infty \top}(t){\bm f}_u(\bm{x}^\infty(t),\bm{u}^\infty (t)) +  {\tilde{\bm{\mu}}}^{\infty \top}(t){\bm h}_u(\bm{x}^\infty(t),\bm{u}^\infty (t)) \right\Vert = 0 \, ,
\end{split}
\end{equation*}
for all $t\in[0,1]$. This proves that ${\bm x}^\infty(t),{\bm u}^\infty(t), \tilde{\bm \lambda}^\infty(t),\tilde{\bm \mu}^\infty(t)$ and $\bar{\bm \nu}^\infty(t)$ satisfy Equations \eqref{eq:Plambda2}-\eqref{eq:Plambda5}.

\qedproof

\begin{thm}[Covector Mapping Theorem] \label{thm:covector}
Under the same assumptions of Theorems \ref{thm:existence_dual} and \ref{thm:consistency_dual}, when $N \to \infty$, the \emph{covector mapping}
$$
\bm{x}_{N}^*(t) \mapsto \bm{x}^*(t) \, , \quad \bm{u}_{N}^*(t) \mapsto \bm{u}^*(t) \, ,
$$
$$
\frac{\bm{\lambda}_{N}^*}{w} \mapsto \bm{\lambda}^*(t) \, , \quad  \frac{\bm{\mu}_{N}^*(t)}{w} \mapsto \bm{\mu}^*(t) \, , \quad \bar{\bm \nu}^* \mapsto \bm{\nu}^* \, 
$$
is a bijective mapping between the solution of Problem $P_{N\lambda}^{\text{clos}}$ and the solution of Problem $P_\lambda$.

\qedlem
\end{thm}

\textbf{Proof:} The above result follows directly from Theorems \ref{thm:existence_dual} and \ref{thm:consistency_dual}.

\qedproof

\section{Numerical examples} \label{sec:numericalresults}
This section presents four numerical examples aimed at validating the convergence properties of the proposed method based on Bernstein approximation. In the first example we consider a non-linear one-dimensional optimal control problem with smooth state, input and costate. In the second example we investigate the applicability and efficacy of our approach, when solving a bang-bang optimal control problem. In the third and fourth examples we demonstrate the benefits of the proposed method, when solving trajectory generation problems for single-vehicle and multi-vehicle missions, respectively. The results are obtained using MATLAB's built in \emph{fmincon} function.

\subsection{Example 1}
\begin{figure}
	\centering
	\includegraphics[width=0.5\textwidth]{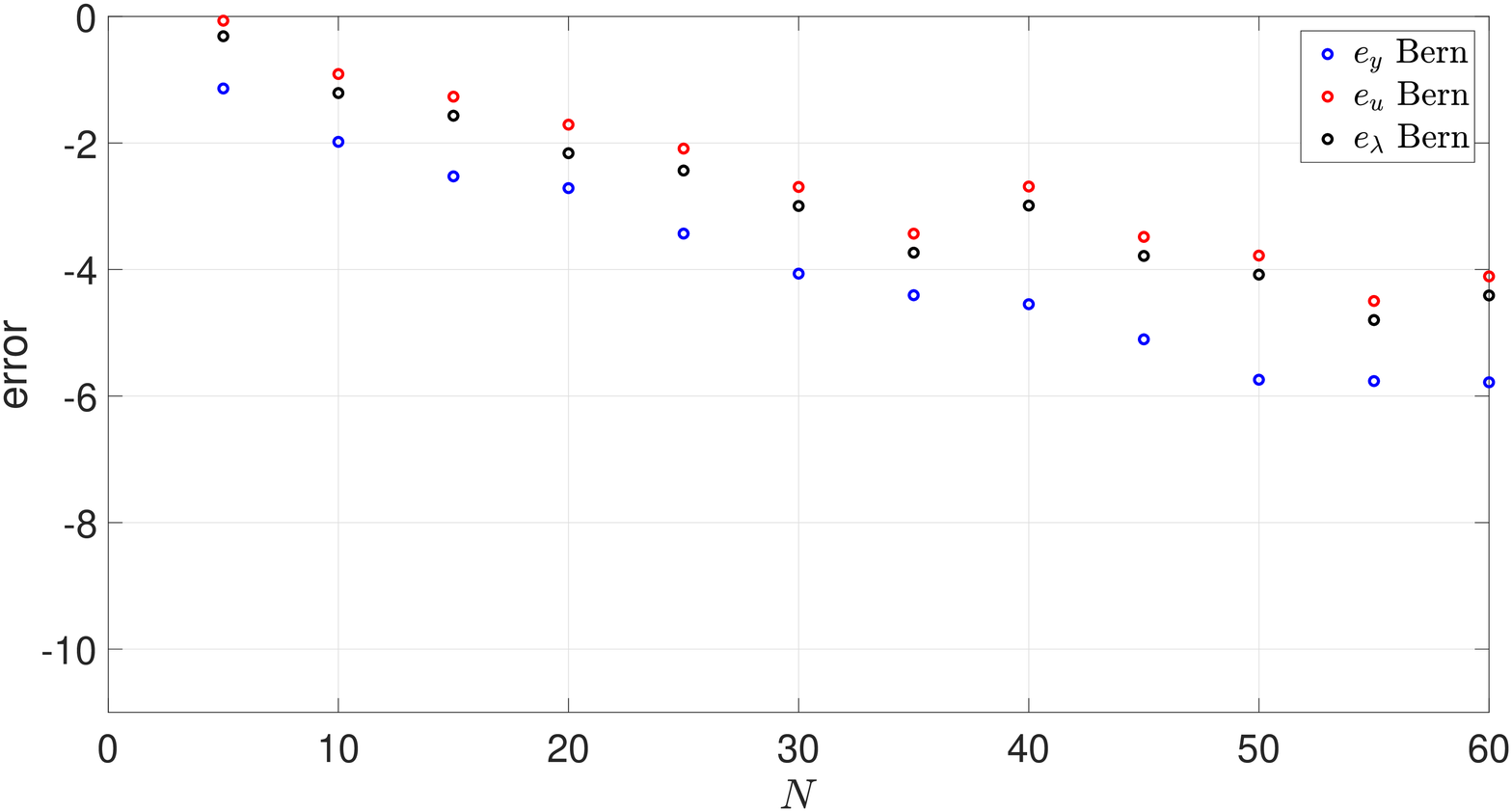}
	\caption{Error in Bernstein approximation method for Example \ref{ex:ex1}.}
	\label{fig:ex1conv_Bern}
\end{figure}

\begin{figure}
	\centering
	\includegraphics[width=0.5\textwidth]{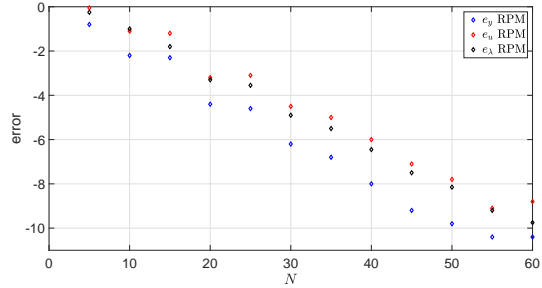}
	\caption{Error in Radau Pseudospectral Method (RPM) for Example \ref{ex:ex1} (taken from \cite{garg2011direct}).}
	\label{fig:ex1conv_PS}
\end{figure}

\begin{figure}
	\centering
	\includegraphics[width=0.5\textwidth]{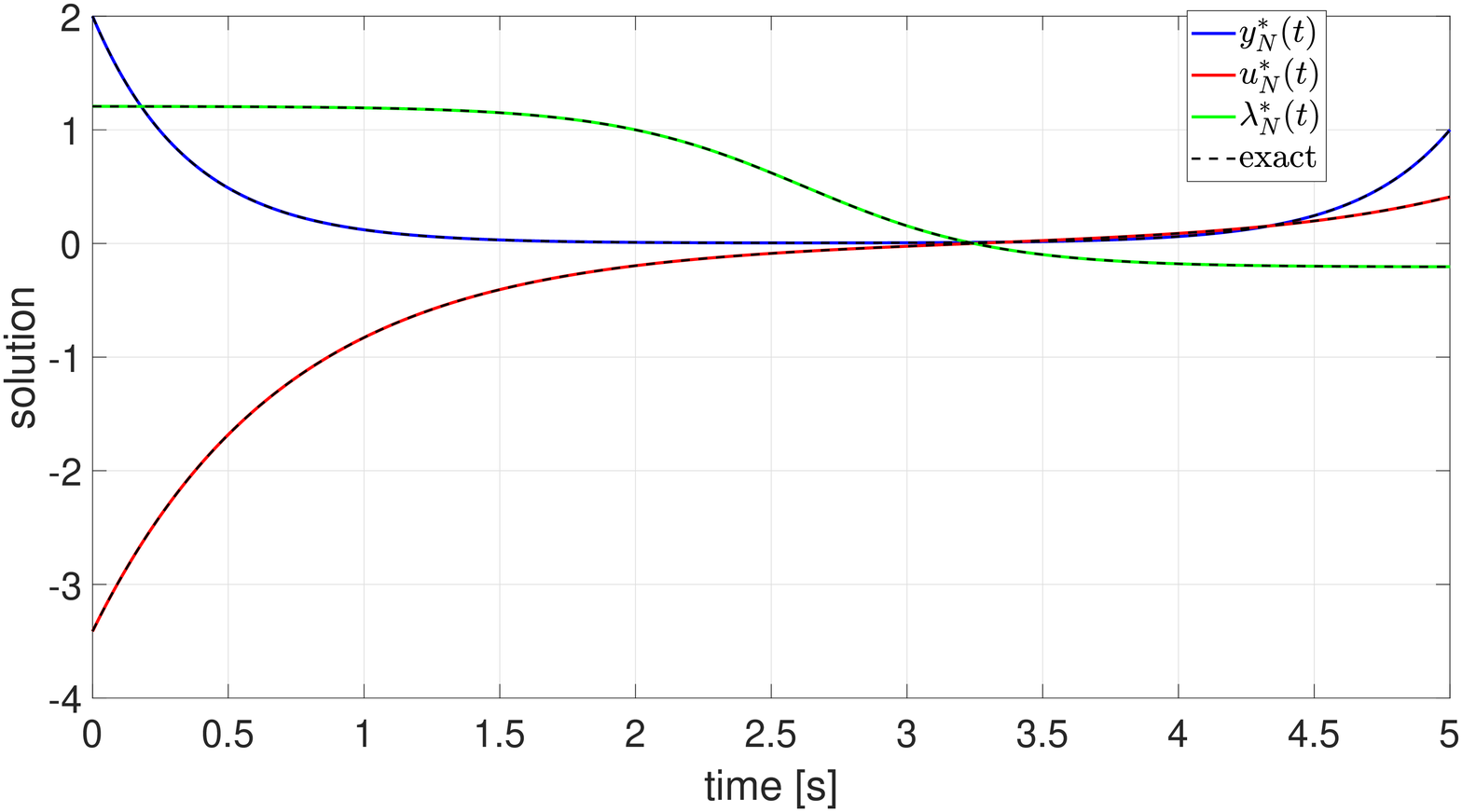}
	\caption{Solution to Example \ref{ex:ex1} using the Bernstein approximation method for order $N=40$.}
	\label{fig:ex1output}
\end{figure}

Consider the following optimal control problem taken from \cite{garg2011direct}:
\begin{exm} \label{ex:ex1}
Determine $y: [0,5] \to \IR$ and $u: [0,5] \to \IR$ that minimize
\begin{equation*} \label{eq:costfunc_ex1}
I(y(t),u(t)) = \frac{1}{2} \int_0^{5} (y(t)+u^2(t))dt \,,
\end{equation*}
subject to
\begin{align*}
& \dot{y}(t) =  2y(t) + 2u(t)\sqrt{y(t)} \, , \quad \forall t \in[0,5], \nonumber \\ 
& y(0) = 2 \, , \nonumber \\ 
& y(5) = 1 \, . \nonumber \\ 
\end{align*}
\end{exm}
The above example was solved using the Bernstein approximation method for orders $N=5,10,\ldots,55,60$.
Similar to \cite{garg2011direct}, we define the following errors
\begin{equation*}
\begin{split}
e_y & = \max_{k = 0,\ldots,N} \log_{10} |y^*_N(t_k) - y^*(t_k))| \, , \\
e_u & = \max_{k = 0,\ldots,N} \log_{10} |u^*_N(t_k) - u^*(t_k))| \, , \\
e_\lambda & = \max_{k = 0,\ldots,N} \log_{10} |\lambda^*_N(t_k) - \lambda^*(t_k))| \, ,
\end{split}
\end{equation*}
where $y^*_N(t_k)$, $u^*_N(t_k)$ and $\lambda^*_N(t_k)$ are the state, input and costate evaluated at the equidistant time nodes $t_k = 5k/N$, $k=0,\ldots,N$, while $y^*(t_k)$, $u^*(t_k)$ and $\lambda^*(t_k)$ are the exact solutions.
Figure~\ref{fig:ex1conv_Bern} illustrates the convergence of the above errors to zero as the order of approximation increases.
As expected, the Bernstein approximation method is outperformed by pseudospectral methods in terms of converge rate. For example, Figure \ref{fig:ex1conv_PS} presents the results obtained using the Radau Pseudospectral Method (RPM) presented in \cite{garg2011direct}. The data depicted in Figure \ref{fig:ex1conv_PS} are reported from \cite{garg2011direct}, where the authors solved Example \ref{ex:ex1} using RPM with the software OptimalPrime and the NLP solver SNOPT.
Finally, Figure \ref{fig:ex1output} depicts the state, input, and costate obtained using the Bernstein approximation method for order $N=40$ alongside the exact solution.


\subsection{Example 2 - Bang-bang control}
\begin{figure*}
	\centering
	\includegraphics[width=0.9\textwidth]{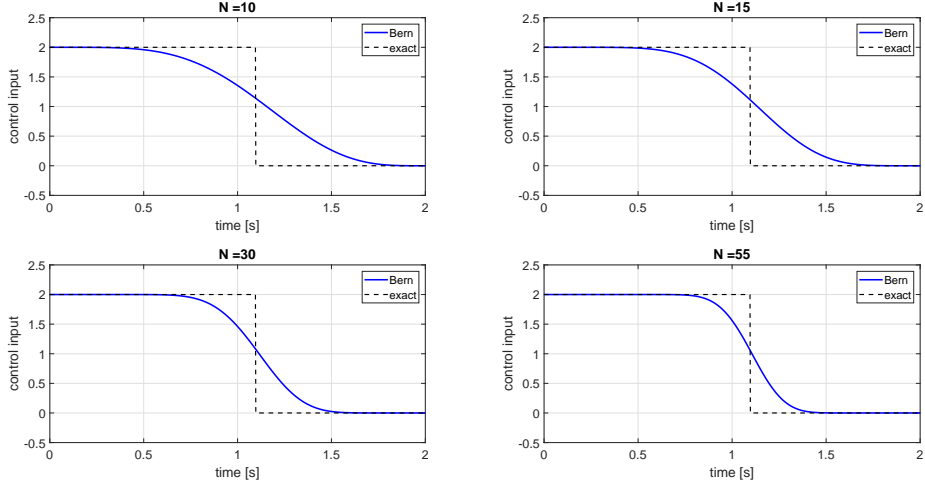}
	\caption{Solution to Example \ref{ex:ex2} using the Bernstein approximation method.}	\label{fig:example2}
\end{figure*}

This example investigates the efficacy of the proposed method when dealing with a bang-bang optimal control problem. Consider the following problem:
\begin{exm} \label{ex:ex2}
Determine $y: [0,2] \to \IR$ and $u: [0,2] \to \IR$ that minimize
\begin{equation} \label{eq:costfunc_ex2}
I(y(t),u(t)) = \int_0^{2} (3u(t)-2y(t))dt \,,
\end{equation}
subject to
\begin{align}
& \dot{y}(t) =  y(t) + u(t) \, , \quad \forall t \in[0,2], \label{eq:dynamicconstraint_ex2} \\
& y(0) = 4 \, , \label{eq:equalityconstraint1_ex2} \\
& y(2) = 39.392 \, , \label{eq:equalityconstraint2_ex2}  \\
& 0 \leq u(t) \leq 2 \quad \forall t \in[0,2]\, .
\end{align}
\end{exm}
The optimal control for the above example is:
$$
u^*(t) =
\begin{cases}
&2 \qquad 0 \leq t \leq 1.096 \\
&0 \qquad 1.096 \leq t \leq 2 .
\end{cases}
$$
Example \ref{ex:ex2} is solved using the Bernstein approximation method for orders of approximation $N=10,15,30,55$. The results are illustrated in Figure \ref{fig:example2}.
It can be noted that the solutions resulting from the Bernstein approximation method converge, albeit slowly, to the exact solution, and behave nicely despite the discontinuity. The solutions have no jumps in the neighborhood of the discontinuities (the Gibbs phenomenon does not occur for function approximation via Bernstein polynomials \cite{gzyl2003approximation}), and the exact value of the discontinuity ($t=1.096 \rm{s}$) is detected with reasonable accuracy even for low orders of approximation. The reader is referred to \cite{tohidi2013efficient}, where the discretization method presented in \cite{tohidi2011legendre} is implemented to solve Example \ref{ex:ex2}, emphasizing the inefficacy of pseudospectral methods when approximating bang-bang solutions. Additional references discussing the performance of pseudospectral methods when dealing with bang-bang optimal control problems can be found in \cite{ross2004pseudospectral,darby2011hp,Wei2005Discontinuous,gelb2006robust}.

\subsection{Example 3 - Trajectory generation for a single vehicle}
In this section, the 2D trajectory generation problem for a single vehicle is considered. The vehicle, modelled as a single integrator,
is required to navigate from the initial position $\bm{x}_0 = [-500, -900]$m to the final destination $\bm{x}_f = [1500, -600]$m, while minimizing the time of arrival. The algorithm must ensure a minimum separation of $E=50$m with three obstacles positioned at $\bm{p}_{o_1} = [ 0 -800]^\top$m, $\bm{p}_{o_2} = [450 -750]^\top$m and $\bm{p}_{o_3} = [850 -730]^\top$m. Finally, the norm of the input must remain within minimum and maximum saturation limits $u_{\min}=15$m/s and $u_{\max}=32$m/s. This problem can be formally stated as follows:

\begin{exm} \label{ex:ex3}
Determine $\bm{x}(t),\bm{u}(t)$ and $t_f$ that minimize
\begin{equation*}
I(\bm{x}(t),\bm{u}(t)) = \int_0^{t_f} dt \,,
\end{equation*}
subject to
\begin{align*}
& \dot{\bm{x}}(t) = \bm{u}(t)\, , \quad \forall t \in[0,t_f],  \\
& \bm{x}(0)=\bm{x}_0 \, , \quad \bm{x}(t_f)=\bm{x}_f, \\
& ||\bm{x}(t)-\bm{p}_{o_i}|| \geq E  \, , \quad \forall t\in [0,t_f], \quad i=1,2,3  \, , \\
& u_{\min} \leq ||\bm{u}(t)|| \leq u_{\max} \, , \quad \forall t\in [0,t_f].
\end{align*}
\end{exm}

\begin{figure*}
	\centering
	\includegraphics[width=0.9\textwidth]{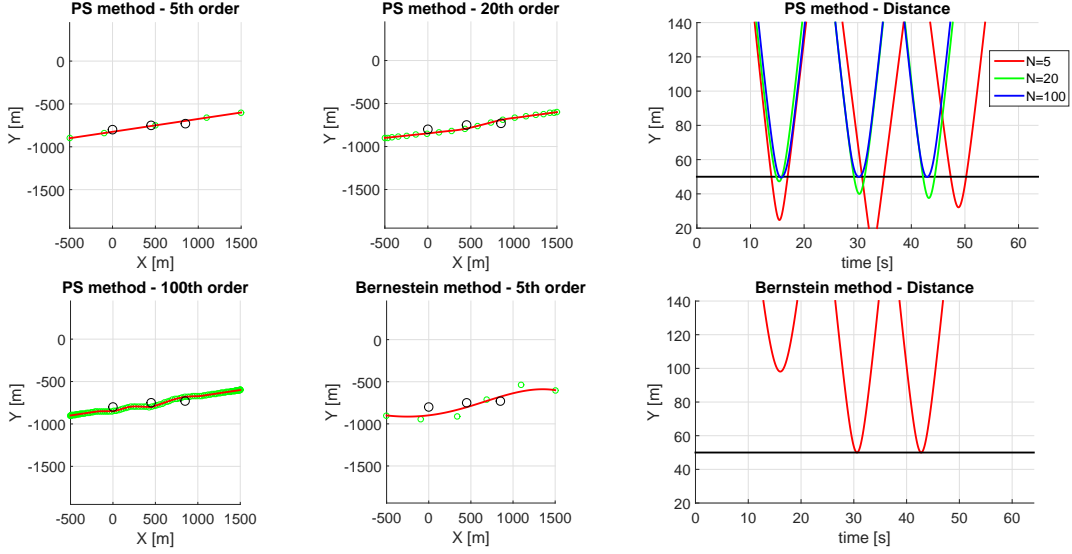}
	\caption{Legendre pseudospectral vs Bernstein approximation method: collision avoidance with multiple obstacles.}
	\label{fig:1veh}
\end{figure*}

The discretization method proposed in this paper is compared to the Legendre pseudospectral method \cite{gong2007pseudospectral}. The results are enclosed in Figure \ref{fig:1veh}.
The top-left, top-center and bottom-left figures show the trajectories, obtained using the pseudospectral method with orders of approximation 5, 20 and 100, respectively. As discussed in the introduction, the pseudospectral method enforces the constraints only at the discretization nodes and not in between them. By increasing the number of nodes, the distance between the entire trajectory and the obstacles increases towards the desired value $E = 50m$. However, as demonstrated by the top-right figure, which depicts the distance between the trajectories and the obstacles for the three order of approximations indicated above, the minimum separation constraint is never satisfied.
On the other hand, the bottom-center figure shows that with the proposed method, even by choosing a small number of nodes (N=5 in this example), the collision avoidance constraint can be computed for the entire curve using, for example, the minimum distance algorithm introduced in Section \ref{sec:MathematicalPreliminaries}, Property \ref{app.prop:mindist}. Thus, collision avoidance is satisfied for the entire trajectory. The bottom-right figure supports this claim by showing that the distance between the trajectory and the obstacles is always greater than the required value.

\subsection{Example 4 - Trajectory generation for multi-vehicle missions}
As evidenced by the previous example, the possibility of choosing low order of approximations while guaranteeing constraint satisfaction is the strength of the proposed approach, which prioritizes safety and feasibility of the trajectories over optimality.
The advantage of the method becomes more evident in multi-vehicle missions, where the trajectories assigned to the vehicles must be (spatially or temporally)\footnote{Spatial separation is guaranteed if the minimum distance between any two points on two paths is greater than or equal to a minimum spatial clearance (the paths never intersect); temporal separation is achieved if for any time $t$ the minimum distance between two vehicles is greater than or equal to the minimum spatial clearance.} separated.
Consider for example a mission scenario in which $n$ vehicles, starting from their initial positions, have to follow spatially separated trajectories to reach predefined final destinations. By adopting the pseudospectral method described above, spatial separation would have to be enforced by imposing separation constraints between every node of every trajectory. Thus, the problem would have $ \left(\begin{matrix} n \\ 2 \end{matrix}\right)  N^2$ separation constraints, where $N$ is the number of nodes, and $\left(\begin{matrix} n \\ 2 \end{matrix}\right)$ is the binomial coefficient. An increased number of nodes (dictated, perhaps, by reasons similar to the ones discussed in the previous example) would increase the complexity in the search for the optimal solution, making the pseudospectral approach practically infeasible for these types of applications.
On the other hand, with the approach of this paper the constraint satisfaction is achieved independently of the number of nodes. This is discussed in the next simulation scenario.

\begin{figure}[h]
	\centering
	\includegraphics[width=0.5\textwidth]{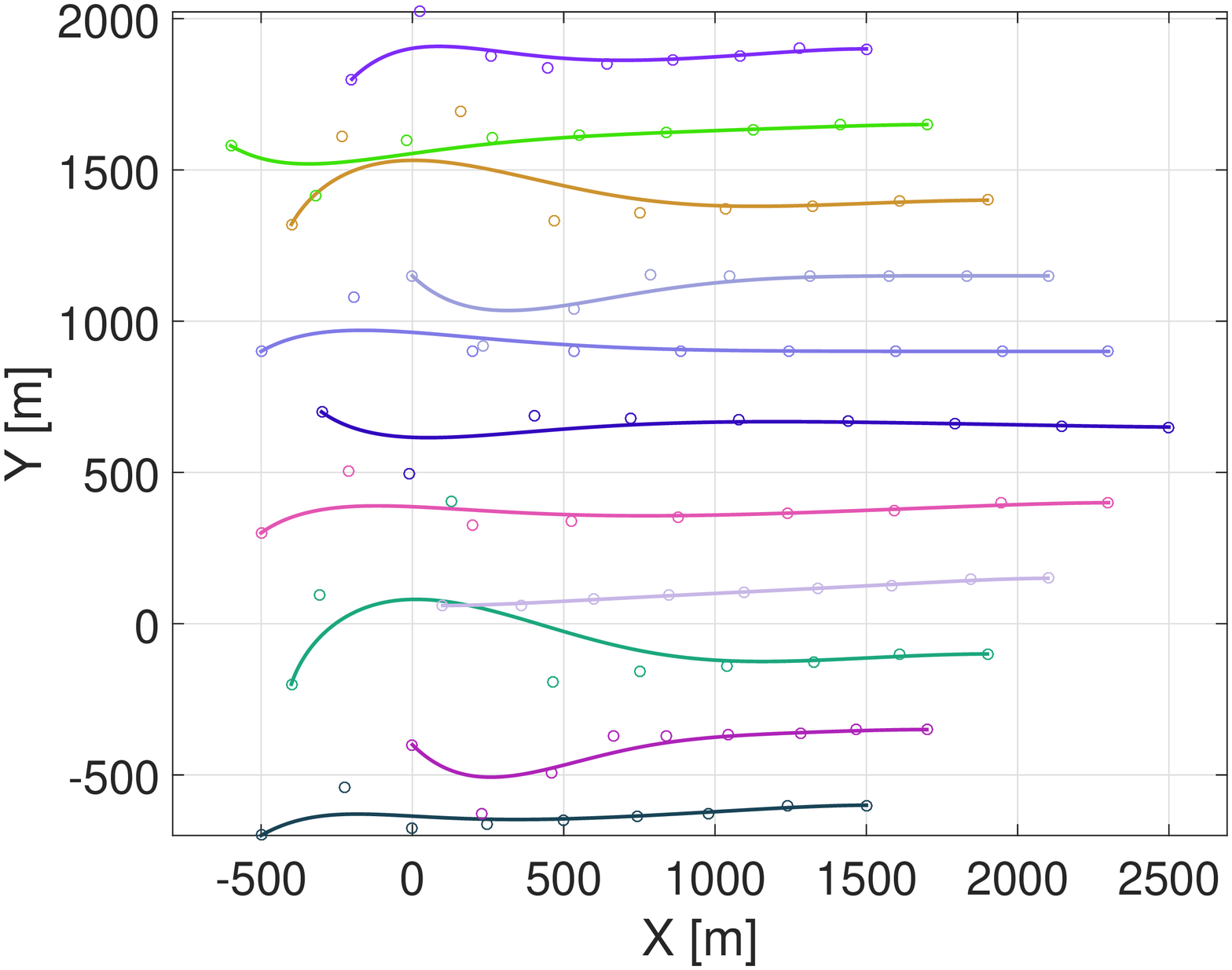}
	\caption{Multiple vehicles mission - temporal separation: 2D paths.}
	\label{fig:11veh_2Dpaths_temp}
\end{figure}

Figure~\ref{fig:11veh_2Dpaths_temp} illustrates the results of a multi-vehicle mission in which $n=11$ vehicles, starting from their initial positions, have to reach a `$>$' shaped formation while minimizing the time of arrival. The Bernstein approximation method is employed with the order of approximation $N=8$. The dynamics of the $i$th vehicle are governed by the following differential equations
\begin{equation}
\begin{cases}
\dot{x}_{1,i}(t) & = V_i(t) \cos(x_{3,i}(t)) \\
\dot{x}_{2,i}(t) & = V_i(t) \sin(x_{3,i}(t)) \\
\dot{x}_{3,i}(t) & = \omega_i(t) \, ,
\end{cases}
\end{equation}
with input $\bm{u}_i(t) = [V_i(t) \, , \, \omega_i(t)]^\top$ and state $\bm{x}_i(t) = [x_{1,i}(t) \, , \, x_{2,i}(t) \, , \, x_{3,i}(t)]^\top$. The input constraints are:
\begin{equation} \label{eq:speedconstraint}
V_{\min}^2 \leq V_i^2(t) \leq V_{\max}^2 \, ,
\end{equation}
\begin{equation} \label{eq:angrateconstraint}
-\omega_{\max} \leq \omega_i(t) \leq \omega_{\max} \, ,
\end{equation}
with $V_{\min} = 15 \rm{m/s}$, $V_{\max} = 32 \rm{m/s}$ and $\omega_{\max} = 0.3 \rm{rad/s}$.
Finally, temporal separation constraints are imposed between each pair of trajectories, i.e.
\begin{equation} \label{eq:vehconstraint}
|| \bm{p}_{i}(t)-\bm{p}_{j}(t) || \geq E \, ,
\end{equation}
$ \forall i,j = 1,\ldots,11, \, i\neq j\, \, , \forall t \in [0,t_f]$, where $E=50\rm{m}$ and $\bm{p}_{i}(t) = [x_{1,i}(t) \, , \, x_{2,i}(t)]^\top$.

\begin{figure}[h]
	\centering
	\includegraphics[width=0.5\textwidth]{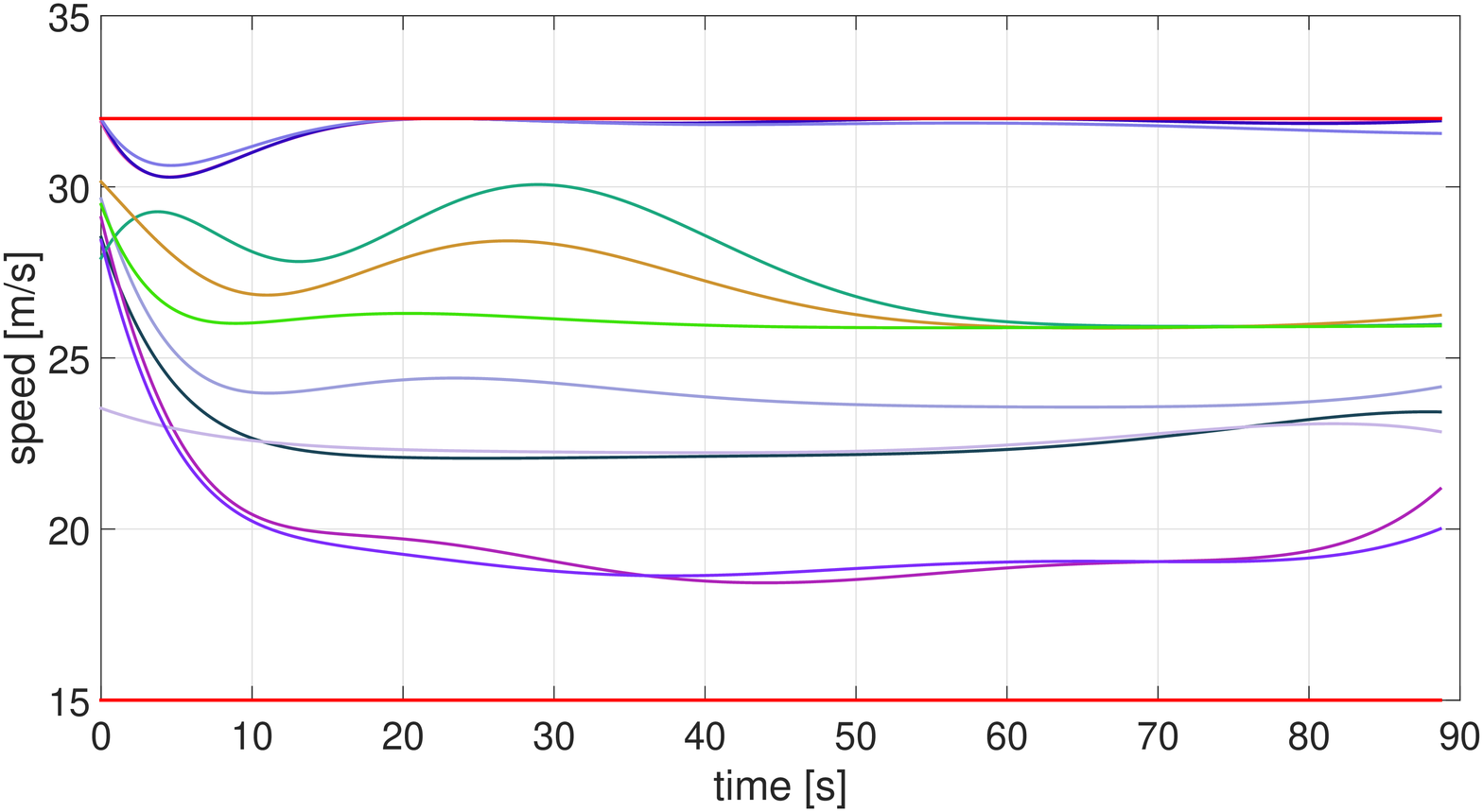}
	\caption{Multiple vehicles mission - temporal separation: speed profiles.}
	\label{fig:11veh_speed_temp}
\end{figure}

\begin{figure}[h]
	\centering
	\includegraphics[width=0.5\textwidth]{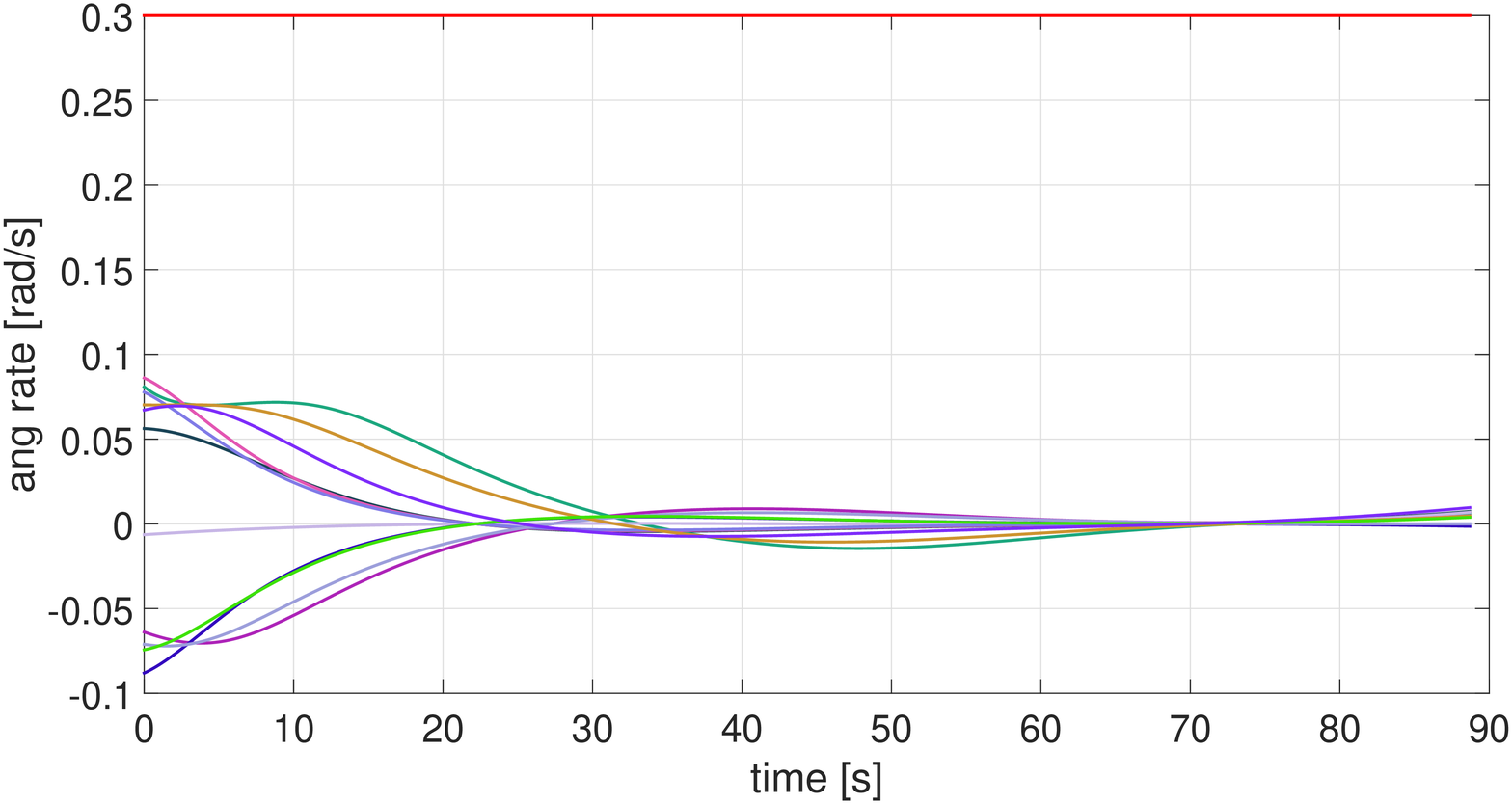}
	\caption{Multiple vehicles mission - temporal separation: angular rates.}
	\label{fig:11veh_angrate_temp}
\end{figure}

The constraints in Equations \eqref{eq:speedconstraint}, \eqref{eq:angrateconstraint} and \eqref{eq:vehconstraint} are computed by exploiting the properties of Bernstein polynomials. In particular, we used the minimum distance and the computation of extrema algorithms from \cite{chang2011computation,ChoePhd}. Figures \ref{fig:11veh_speed_temp} and \ref{fig:11veh_angrate_temp} show the time history of the speeds and angular rates, respectively, demonstrating that the input saturation constraints are satisfied for all times.

\begin{figure}[h]
	\centering
	\includegraphics[width=0.5\textwidth]{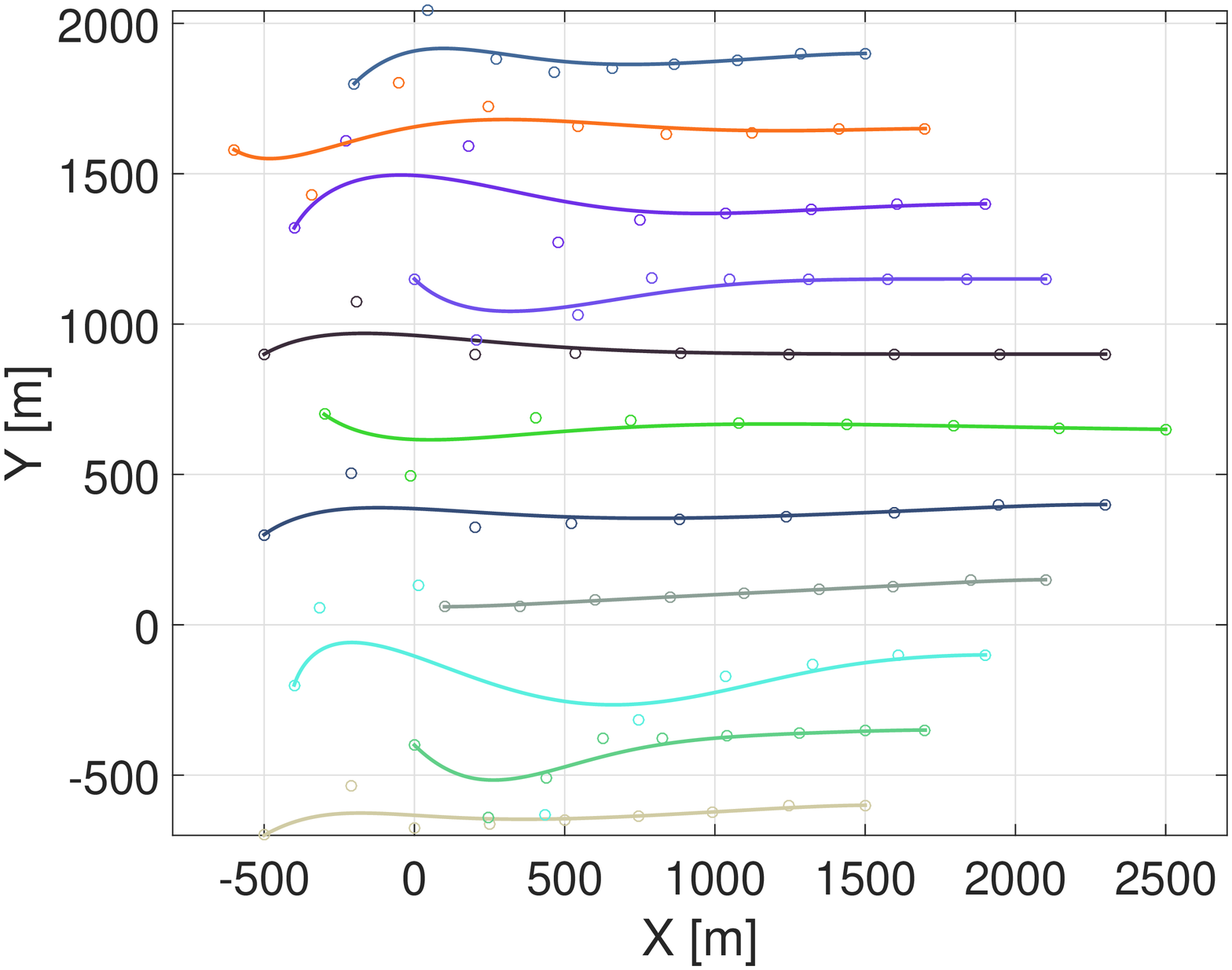}
	\caption{Multiple vehicles mission - spatial separation: 2D paths.}
	\label{fig:11veh_2Dpaths_spat}
\end{figure}

\begin{figure}[h]
	\centering
	\includegraphics[width=0.5\textwidth]{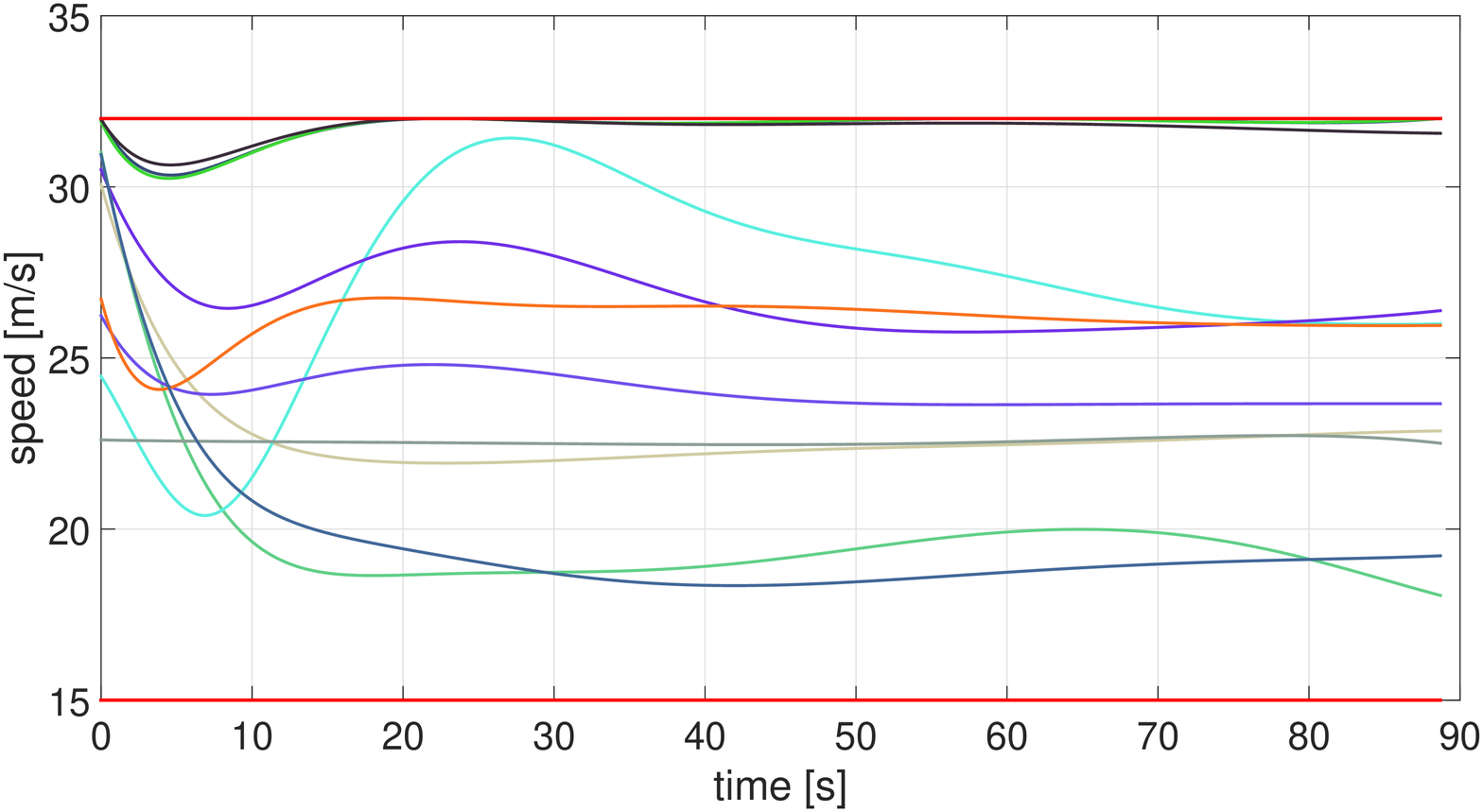}
	\caption{Multiple vehicles mission - spatial separation: speed profiles.}
	\label{fig:11veh_speed_spat}
\end{figure}

\begin{figure}[h]
	\centering
	\includegraphics[width=0.5\textwidth]{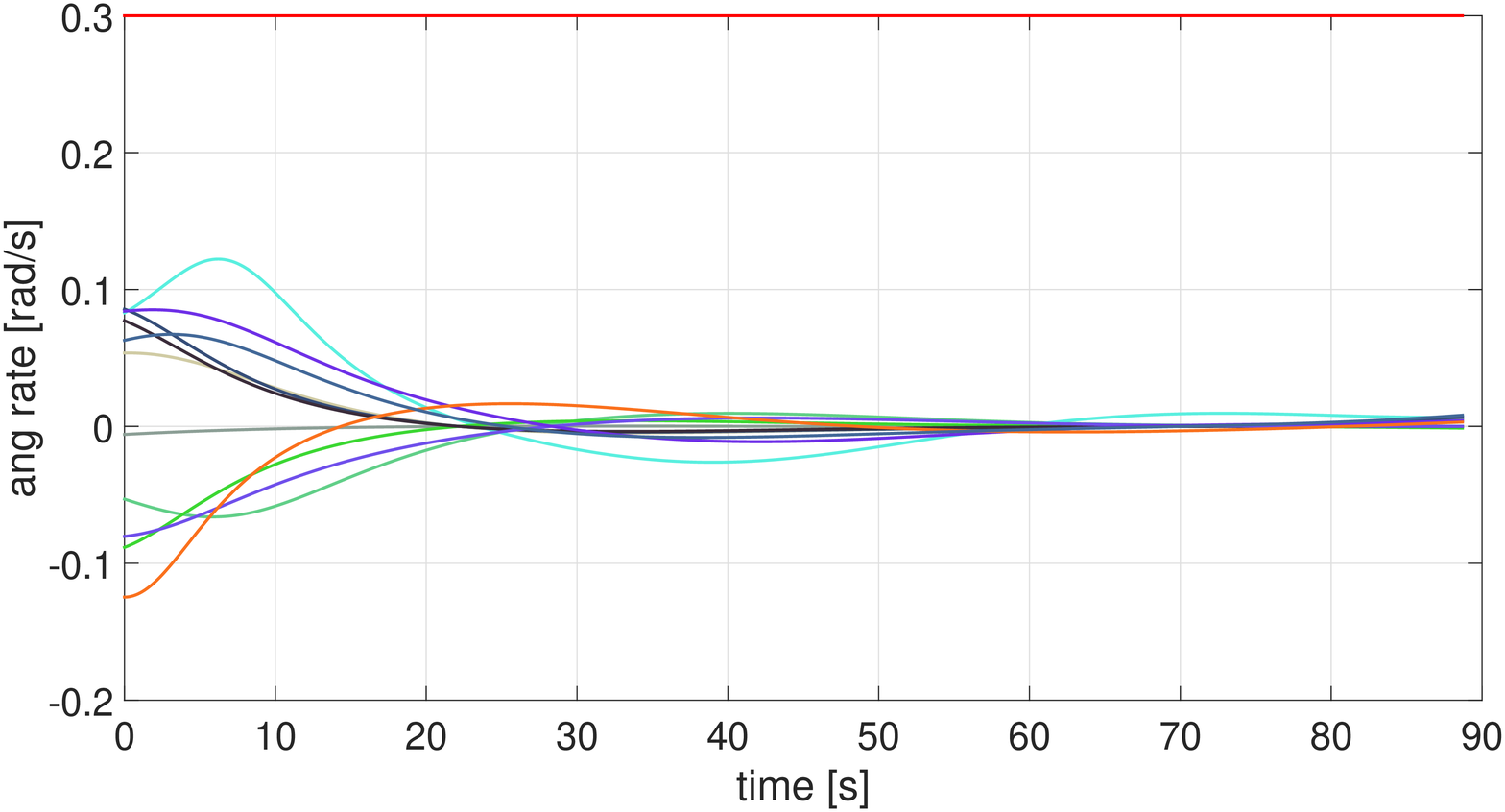}
	\caption{Multiple vehicles mission - spatial separation: angular rates.}
	\label{fig:11veh_angrate_spat}
\end{figure}

At last, the same simulation is repeated, but the temporal separation constraint given by Equation \eqref{eq:vehconstraint} is replaced by the (more stringent) spatial separation requirement
\begin{equation} \label{eq:vehconstraint_spat}
|| \bm{p}_{i}(t_k)-\bm{p}_{j}(t_p) || \geq E \, ,
\end{equation}
$ \forall i,j = 1,\ldots,11, \, i\neq j\, \, \forall t_k,t_p \in [0,t_f]$.
Figure \ref{fig:11veh_2Dpaths_spat} depicts the 2D trajectories. Figures \ref{fig:11veh_speed_spat} and \ref{fig:11veh_angrate_spat} illustrate the speed and angular rate commands, respectively.



\section{Conclusions} \label{sec:conclusions}
This paper proposed a numerical method to approximate nonlinear constrained optimal control problems by nonlinear programming (NLPs) using Bernstein polynomial approximation.
A rigorous analysis is provided that shows convergence of the solution of the NLP to the solution of the continuous-time problem. A set of conditions are derived under which the Karush-Kuhn-Tucker multipliers of the NLP converge to the costates of the optimal control problem. This led to the formulation of the Covector Mapping Theorem for Bernstein approximation, which enables numerical computation of the costates.
The theoretical findings are validated through several numerical examples, and the advantages and disadvantages of the proposed method are discussed.


\appendix
\section*{Appendix}
\small

\section{Proof of Equation \eqref{eq:inequalities1}}
\label{app:inequalities1}
Let us focus on Equation \eqref{eq:inequalities11}. Adding and subtracting $\int_0^1 F_x(\bm{x}_N(t),\bm{u}_N(t)) b_{0,N}(t) dt$, we have
\begin{equation} \label{eq:appeqn1}
\begin{split}
& \left\Vert w \sum_{j = 0}^N F_x(\bm{x}_N(t_j),\bm{u}_N (t_j)) b_{0,N}(t_j) - \int_0^1 F_x(\bm{x}(t),\bm{u}(t)) b_{0,N}(t) dt \right\Vert \leq \\
& \qquad \left\Vert w \sum_{j = 0}^N F_x(\bm{x}_N(t_j),\bm{u}_N (t_j)) b_{0,N}(t_j) - \int_0^1 F_x(\bm{x}_N(t),\bm{u}_N(t)) b_{0,N}(t) dt + \int_0^1 F_x(\bm{x}_N(t),\bm{u}_N(t)) b_{0,N}(t) dt \right. \\
& \qquad \quad \left. - \int_0^1 F_x(\bm{x}(t),\bm{u}(t)) b_{0,N}(t) dt \right\Vert \\
& \quad \leq  \left\Vert w \sum_{j = 0}^N F_x(\bm{x}_N(t_j),\bm{u}_N (t_j)) b_{0,N}(t_j) - \int_0^1 F_x(\bm{x}_N(t),\bm{u}_N(t)) b_{0,N}(t) dt \right\Vert  \\
& \qquad \quad + \left\Vert \int_0^1 F_x(\bm{x}_N(t),\bm{u}_N(t)) b_{0,N}(t) dt - \int_0^1 F_x(\bm{x}(t),\bm{u}(t)) b_{0,N}(t) dt \right\Vert.
\end{split}
\end{equation}
Using Lemma \ref{lem:quadrature} and continuity of $F_x(\bm{x}_N(t),\bm{u}_N (t))$ and $b_{0,N}(t)$, the first term on the right hand side of the inequality above satisfies
$$
\left\Vert w \sum_{j = 0}^N F_x(\bm{x}_N(t_j),\bm{u}_N (t_j)) b_{0,N}(t_j) - \int_0^1 F_x(\bm{x}_N(t),\bm{u}_N(t)) b_{0,N}(t) dt \right\Vert \leq C_I W_{F_xb_{0,N}}(N^{-\frac{1}{2}}) ,
$$
where $W_{F_xb_{0,N}}(\cdot)$ is used to denote the modulus of continuity of the product $F_x(\bm{x}_N(t),\bm{u}_N(t)) b_{0,N}(t)$, with $F_x(\bm{x}_N(t),\bm{u}_N(t))$ being a bounded function due to its continuity over a bounded domain. Denote its bound as $F_{x,\max}$. Notice that also $b_{0,N}(t)$ is bounded, as $\max_{t\in[0,1]}b_{0,N}(t) \leq 1$. Then using the properties of the modulus of continuity, we get
\begin{equation} \label{eq:appeqn11}
\begin{split}
\left\Vert w \sum_{j = 0}^N F_x(\bm{x}_N(t_j),\bm{u}_N (t_j)) b_{0,N}(t_j) - \int_0^1 F_x(\bm{x}_N(t),\bm{u}_N(t)) b_{0,N}(t) dt \right\Vert & \leq C_I F_{x,\max}W_{b_{0,N}}(N^{-\frac{1}{2}}) + C_I W_{F_x}(N^{-\frac{1}{2}}) \\
& \leq C_I F_{x,\max}N^{-\frac{1}{2}} + C_I W_{F_x}(N^{-\frac{1}{2}}) ,
\end{split}
\end{equation}
where $W_{F_x}(\cdot)$ is the modulus of continuity of $F_x$, and $C_I$ is a positive constant independent of $N$. Furthermore, we have
\begin{multline} \label{eq:appeqn12}
\left\Vert \int_0^1 F_x(\bm{x}_N(t),\bm{u}_N(t)) b_{0,N}(t) dt - \int_0^1 F_x(\bm{x}(t),\bm{u}(t)) b_{0,N}(t) dt \right\Vert \\ \leq \int_0^1 \left\Vert F_x(\bm{x}_N(t),\bm{u}_N(t)) b_{0,N}(t)  - F_x(\bm{x}(t),\bm{u}(t)) b_{0,N}(t) \right\Vert dt  \leq L_{F_x} (C_x W_x(N^{-\frac{1}{2}})+C_u W_u(N^{-\frac{1}{2}})) \, ,
\end{multline}
where $L_{F_x}$ is the Lipschitz constant of $F_x$, $C_x < 5n_x/4$, $C_u < 5n_u/4$, and $W_x(\cdot)$ and $W_u(\cdot)$ are the moduli of continuity of $\bm{x}$ and $\bm{u}$, respectively.
Combining Equations \eqref{eq:appeqn11} and \eqref{eq:appeqn12} with Equation \eqref{eq:appeqn1}, yields
\begin{multline}
\left\Vert w \sum_{j = 0}^N F_x(\bm{x}_N(t_j),\bm{u}_N (t_j)) b_{0,N}(t_j) - \int_0^1 F_x(\bm{x}(t),\bm{u}(t)) b_{0,N}(t) dt \right\Vert \\ \leq  C_I F_{x,\max}N^{-\frac{1}{2}} + C_I W_{F_x}(N^{-\frac{1}{2}}) + L_{F_x} (C_x W_x(N^{-\frac{1}{2}})+C_u W_u(N^{-\frac{1}{2}})) \, ,
\end{multline}
which proves the bound in Equation \eqref{eq:inequalities11}. The bounds in Equations \eqref{eq:inequalities12}-\eqref{eq:inequalities14} follow easily using an identical argument.


\bibliographystyle{IEEEtran}
\bibliography{IEEEtran}

\end{document}